\documentclass{amsart}
\usepackage{amsfonts,amssymb,amsmath,a4wide}
\newcommand{\qedbox}{ \fbox{}}

\newtheorem{teo}{Theorem}[section]
\newtheorem{lema}{Lemma}[section]
\newtheorem{pro}{Proposition}[section]
\newtheorem{defi}{Definition}[section]\newtheorem{rema}{Remark}[section]
\newtheorem{coro}{Corollary}[section]
\numberwithin{equation}{section}\def\ndt{\noindent}

\def\A{\mathcal{A}}
\def \V{\mathcal{V}}
\def \H{\mathcal{H}}

\def\og{\overline{\gamma}}
\def \bnabla{\overline{\nabla}}

\def \z{\zeta}

\def\J{\mathcal{J}}
\def\w{\wedge}
\def \l4{[[\Lambda^{0,4}]]}
\def \ll3{[[\Lambda^{0,3}]]}
\def \ho{\otimes_1}
\def \ko{\otimes_2}
\def \bJ{\mathbb{J}}

\def\Ric{\text{Ric}}
\def\dim{\text{dim\ }}
\def\sideremark#1{\ifvmode\leqslantavevmode\fi\vadjust{\vbox to0pt{\vss
 \hbox to 0pt{\hskip\hsize\hskip1em
 \vbox{\hsize2.5cm\tiny\raggedright\pretolerance10000
 \noindent #1\hfill}\hss}\vbox to8pt{\vfil}\vss}}}%

\begin{document}\larger[2]
\title[]{On the cohomology algebra of some classes of geometrically formal manifolds}
\author[J.-F. Grosjean]{J.-F. Grosjean}
\author[P.- A. Nagy]{P.-A. Nagy}
\address[J.-F. Grosjean]{Institut \'Elie Cartan, Universit\'e H. Poincar\'e,
Nancy I, B.P.239 F-54506 Vandoeuvre-L\`es-Nancy, France}
\email{grosjean@iecn.u-nancy.fr}
\address[P.-A. Nagy]{Department of Mathematics, University of Auckland, Private Bag 92019, Auckland, New Zealand}
\email{nagy@math.auckland.ac.nz}
\date{\today}

\begin{abstract} We investigate harmonic forms of geometrically formal metrics, which are defined as those having the 
exterior product of any two harmonic forms still harmonic. 
We prove that a formal Sasakian metric can exist only on a real cohomology sphere and that holomorphic forms of a formal 
K\"ahler metric are parallel w.r.t. the Levi-Civita connection. In the general Riemannian case a formal metric with maximal second Betti number is 
shown to be flat . Finally we prove that a six-dimensional manifold with $b_1 \neq 1, b_2 \geqslant 2$ and not having the real cohomology algebra of $\mathbb{T}^3 \times S^3$ carries a symplectic structure as soon as it admits a formal metric.
\end{abstract}
\maketitle
\tableofcontents
\section{Introduction}
Let $(M^n,g)$ be a compact oriented Riemannian manifold. We denote by $\Lambda^pM, 0 \leqslant p \leqslant n$ the space 
of smooth, real-valued, differential $p$-forms of $M$. We have then a differential complex
$$ \ldots \rightarrow \Lambda^{p}M \stackrel{d}{\rightarrow}  \Lambda^{p+1}M \ldots $$
where $d$ is the exterior derivative. The $p$-th cohomology group of this complex, known as the $p$-th deRham cohomology group will be denoted by $H_{DR}^{p}(M)$. The Riemannian 
metric $g$ induces a scalar product at the level of differential forms, hence one can consider also the operator $d^{\star}$, the formal adjoint of $d$. For $0 \leqslant p \leqslant n$ we define 
the space of harmonic $p$-forms by setting 
$$ {\mathcal{H}}^p(M,g)=\{ \alpha \in \Lambda^{p}M : \Delta \alpha=0 \}.$$
Here the Laplacian $\Delta$ is defined by
$$ \Delta=dd^{\star}+d^{\star}d.$$
Classical Hodge theory produces an isomorphism
\begin{equation} \label{hodge} H_{DR}^{p}(M) \cong {\mathcal{H}}^{p}(M,g)
\end{equation}
for all $0 \leqslant p \leqslant n$. Whilst $H^{\star}(M)=\displaystyle \bigoplus_{p\geqslant 0} H^p_{DR}(M)$ is a graded algebra, generally 
${\mathcal{H}}^{\star}=\displaystyle\bigoplus_{p\geqslant 0}{\mathcal{H}}^{p}(M,g)$ is not an algebra with respect to the wedge product operation for there is no reason that the 
isomorphism \ref{hodge} descends to the level of harmonic forms. Our next definition is related to this fact.
\begin{defi} \cite{Kot}Let $M^n$ be a compact and oriented manifold. 
\begin{itemize}
\item [(i)]A Riemannian metric $g$ on $M$ is formal if the exterior product of any two harmonic (w.r.t. $g$) forms remains harmonic;
\item[(ii)] $M$ is geometrically formal if it admits a formal metric.
\end{itemize}
\end{defi}
A closely related notion is that of topological formality (see \cite{deligne} for instance), which implies that the rational homotopy type of the manifold is a formal consequence of its cohomology ring \cite{Sul}. From the existence of a formal metric it follows that the underlying manifold 
is topologically formal, and this provides obstructions to the existence to such metrics; for instance they cannot exist on nilmanifolds since those 
have non-trivial Massey products, a fact which is in itself an obstruction to formality \cite{deligne,tanre}. On the other hand, simply connected, compact manifolds of dimension not exceeding $6$ are 
topologically formal \cite{fer1, ill}.
$\\$
Now the existence of formal metrics is more directly related to the geometry of the ambient manifold and known obstructions are related to the length of harmonic forms.
\begin{teo} \label{thK}\cite{Kot}Let $(M^n,g)$ be compact and oriented such that $g$ is a formal metric. Then
\begin{itemize}
\item[(i)]the inner product of any two harmonic forms is a constant function; 
\item[(ii)] $b_p(M) \leqslant \begin{pmatrix} n\\ p\end{pmatrix}$ for all $1 \leqslant p \leqslant n$;
\item[(iii)] if in (ii) equality occurs for $p=1$ then $g$ is a flat metric.
\end{itemize}
\end{teo}
Standard examples of formal metrics are provided by compact symmetric spaces for in this case all harmonic forms must be parallel with respect to the Levi-Civita connection. D. Kotschik 
proved that in dimension $4$ every geometrically formal manifold has the real cohomology algebra of a compact symmetric space. One of the current questions related to the 
notion of geometric formality is then to examine up to what extent this is true in general.

In the context of Sasakian geometry, the odd dimensional analogue of K\"ahler geometry we prove 
\begin{teo} \label{Sasformal}Let $(M^{2n+1},g)$ be a compact Sasakian manifold. If $g$ is a formal metric then $M$ is a real cohomology sphere.
\end{teo}
Next we obtain obstructions to the existence of formal K\"ahler metrics, through the study of their holomorphic forms. In this context topological 
formality is no longer restrictive since any K\"ahler manifold is known to have this property \cite{deligne}. 
\begin{teo}\label{holo} Let $(M^{2n},g,J)$ be a compact K\"ahler manifold such that the metric $g$ is formal. Then every harmonic form $\Omega$ of real type 
$(p,0)+(0,p)$ (hence every holomorphic $p$-form ) is parallel with respect to the Levi-Civita connection. Moreover $\Omega$ 
induces in a canonical way a local splitting of $M$ as the Riemannian product of two compact K\"ahler manifolds $M_1$ and $M_2$ so 
that $\Omega$ is zero on $M_1$, non-degenerate on $M_2$ which is Ricci flat. \end{teo}
\begin{rema}
(i) Theorem \ref{holo} was already proved in \cite{nagy} for $p=2$, using arguments relying heavily on the algebraic structure of the space of 
harmonic $2$-forms. For higher degree forms, such results are no longer available. \\
(ii) If in Theorem \ref{holo} we furthermore assume the metric being locally irreducible and not symmetric, it follows from Berger`s holonomy classification 
theorems (see \cite{SMS}) that 
the only cases when we can have a non-vanishing holomorphic form are when $Hol(g)=Sp(m) (n=2m)$ or $Hol(g)=SU(n)$. \\
(iii)  From the above it also follows that if $M$ admits a locally irreducible K\"ahler and formal metric which is not Ricci flat then the Todd genus satisfies $Td(M)=1$. 
\end{rema}
In the second part of the paper we study general properties of $2$-forms which are harmonic w.r.t. a formal metric. We observe that any such $2$-form diagonalises with 
constant eigenvalues and constant rank eigendistributions. This is extending results from \cite{nagy} to the general Riemannian case and can also be used 
as a starting point to give sufficient conditions, essentially phrased in terms of Betti numbers lower bounds, for a formal 
metric to admit a compatible symplectic form in dimension $6$. We prove
\begin{teo}\label{symp} Let $M^6$ be geometrically formal. If $b_1(M) \neq 1$ and $b_2(M) \geqslant 2$ and moreover $M$ has not the 
real cohomology algebra of $\mathbb{T}^3 \times S^3$ then any formal metric on $M$ admits a compatible symplectic form.
\end{teo}
The above result essentially says that in dimension $6$ a geometrically formal manifold $M$ always carries a symplectic structure compatible with the formal metric with the exception of the cases when 
$b_1(M)=1$ or $b_1(M)\neq 1, b_2(M)=0,1$ or when the real cohomology algebra is that of $\mathbb{T}^3 \times S^3$. This suggests that symplectic techniques 
could be used to investigate, under these conditions, the topology and geometry of these manifolds. In dimension $4$, the existence of symplectic forms on geometrically formal manifolds has been extensively treated in \cite{Kot}.

When $b_2(M) \geqslant 3$ Theorem \ref{symp} follows essentially by algebraic arguments mainly using the above mentioned fact on the diagonalization of harmonic $2$-forms of 
a formal metric. To prove it when $b_2(M)=2$ we first show that the absence 
of a compatible symplectic form forces 
the presence of enough harmonic $3$-forms (actually $b_3(M)=6$ in this case). Then we need to perform a rather delicate local analysis, involving the internal symmetries 
of the set harmonic $3$-forms in order to arrive at $b_1(M) \geqslant 2$, a case which can be ruled out algebraically.\\

In the final part of the paper we are concerned with giving a characterisation of geometrically formal Riemannian manifold with maximal second Betti number. We prove
\begin{teo} \label{max2}Let $M^n$ be geometrically formal with $n\geqslant 3$. If $b_2(M)$ is maximal, that is 
$b_2(M)=\begin{pmatrix} n\\ 2\end{pmatrix}$, then any formal metric on $M$ is flat. 
\end{teo}
This clarifies the equality case in Theorem \ref{thK}, (iii) for degree $2$-forms. Note that the assertion in Theorem \ref{max2} is straightforward when $n$ is odd for if 
$n=2k+1$ the formality and the maximality of $b_2$ imply that $b_{2k}(M)$ is maximal. Hodge duality implies then the maximality of $b_1(M)$ and hence the flatness of the metric (see Section 5 for more details). When 
$n$ is even, our point of departure consists in observing that the metric must admit  a compatible almost K\"ahler structure and then work out this situation within the same 
circle of arguments which have led to the proof of Theorem \ref{holo}. 

To conclude, it would be interesting to have results similar to Theorem \ref{symp} in arbitrary even dimensions and of course to give necessary but also sufficient conditions 
for a geometrically formal metric to admit a compatible symplectic structure. In doing so, the difficulties one faces are related to understanding, at the algebraic level, 
the constraints imposed by geometric formality on forms of degree $\geqslant 3$.
\section{Some algebraic facts} 
Let $(V^{2n},g,J)$ be a Hermitian vector space and let $\Lambda^{\star}V$ be its exterior algebra over the reals. Consider the operator $\J : \Lambda^pV \to \Lambda^pV$ acting on 
a $p$-form $\alpha$ by
\begin{equation*}(\J \alpha)(v_1, \ldots, v_p)=\sum \limits_{k=1}^{p}\alpha(v_1, \ldots, Jv_k, \ldots ,v_p)
\end{equation*}
for all $v_1, \ldots v_p$ in $V$. $\J$ acts as a derivation on $\Lambda^{\star}$ and gives the complex bi-grading of the exterior algebra in the following sense. Let $\lambda^{p,q}V$ be given as the $-(p-q)^2$-eigenspace of $\J^2$. Then 
\begin{equation*}\Lambda^sV=\sum \limits_{p+q=s}^{} \lambda^{p,q}V
\end{equation*}
is an orthogonal, direct sum.  Note that $\lambda^{p,q}V=\lambda^{q,p}V$. Of special importance in our discussion are the spaces $\lambda^{p}V=\lambda^{p,0}V$;  forms $\alpha$ in $\lambda^p$ are 
such that $(X_1, \ldots, X_p) \to \alpha(JX_1, X_2, \ldots, X_p)$ is still an alternating form which equals $p^{-1}\J\alpha$. We shall also use the extension of $J$ to $\Lambda^{\star}V$ given by 
\begin{equation*}
(J\alpha)(v_1, \ldots, v_p)=\alpha(Jv_1, \ldots, Jv_p)
\end{equation*}for all $\alpha$ in $\Lambda^pV$ and $v_1, \ldots, v_p$ in $V$. Let $\lambda^pV \ho \lambda^qV$ be the space of tensors $Q\ : \lambda^pV \to \lambda^qV$ which satisfy 
$$[(\bJ Q)(X_1,\ldots,X_p)](Y_1,\ldots,Y_q)=-[\bJ (Q(X_1,\ldots,X_p))](Y_1,\ldots,Y_q)$$ 
(here $\bJ$ as a map of $\lambda^pV$ stands in fact for 
$p^{-1}\J$).  We also define $\lambda^pV \ko \lambda^qV$ to be the space of tensors $Q : \lambda^pV \to \lambda^qV$ such that $Q \bJ=\bJ Q$.
\begin{lema}\label{l31} Let $a : \lambda^pV \otimes \lambda^qV \to \Lambda^{p+q}V$ be the total antisymmetrisation map. Then:  
\begin{itemize}
\item[(i)]
The image of the restriction of a to $ \lambda^pV \ho \lambda^qV \to \Lambda^{p+q}V$ is contained in $\lambda^{p,q}V$;
\item[(ii)]
The image of the restriction of a to $ \lambda^pV \ko \lambda^qV \to \Lambda^{p+q}V$ is contained in $\lambda^{p+q}V$.
\end{itemize}
\end{lema}
\begin{proof}We shall provide a direct proof, but only for (i), that of (ii) being similar. Pick $Q$ in $ \lambda^pV \ho \lambda^qV$. Then 
\begin{equation*}
a(Q)=\sum \limits_{I=(i_1, \ldots ,i_p)}^{} e_{i_1}^{\flat } \wedge \ldots \wedge e_{i_p}^{\flat} \wedge Q(e_{i_1}, \ldots, e_{i_p}) 
\end{equation*}
where for $v$ in $V$ we denote by $v^{\flat}$ the dual, w.r.t to the metric, $1$-form.  Then 
\begin{equation*}
\begin{split}\J (a(Q))&=\sum \limits_{I=(i_1, \ldots, i_p)}^{} \J(e_{i_1}^{\flat} \wedge \ldots \wedge e_{i_p}^{\flat} ) \wedge Q(e_{i_1}, \ldots, e_{i_p}) \\
&+\sum \limits_{I=(i_1, \ldots ,i_p)}^{} e_{i_1}^{\flat} \wedge \ldots \wedge e_{i_p}^{\flat} \wedge \J Q(e_{i_1}, \ldots, e_{i_p}) .
\end{split}
\end{equation*}
For any $1 \leqslant r \leqslant p$ we compute 
\begin{equation*}
\begin{split}
\sum \limits_{I=(i_1, \ldots ,i_p)}^{} e_{i_1}^{\flat } \wedge \ldots \wedge Je_{i_r}^{\flat} \wedge\ldots \wedge e_{i_p}^{\flat} \wedge Q(e_{i_1}, \ldots , e_{i_p})\\
=-\sum \limits_{I=(i_1, \ldots ,i_p)}^{} e_{i_1}^{\flat } \wedge \ldots \wedge(Je_{i_r})^{\flat}\wedge \ldots \wedge e_{i_p}^{\flat} \wedge Q(e_{i_1}, \ldots , e_{i_p})\\
=\sum \limits_{I=(i_1, \ldots ,i_p)}^{} e_{i_1}^{\flat } \wedge \ldots \wedge e_{i_r}^{\flat} \wedge\ldots \wedge e_{i_p}^{\flat} \wedge Q(e_{i_1}, \ldots, Je_{i_r}, \ldots , e_{i_p})\\
=\sum \limits_{I=(i_1, \ldots ,i_p)}^{} e_{i_1}^{\flat} \wedge \ldots \wedge e_{i_p}^{\flat} \wedge (\bJ Q)(e_{i_1}, \ldots , e_{i_p}).
\end{split}
\end{equation*}On the other side we have $\J Q(e_{i_1}, \ldots , e_{i_p})=q\bJ[Q(e_{i_1}, \ldots, e_{i_p})]=-q (\bJ Q)(e_{i_1}, \ldots , e_{i_p})$ and 
putting all these together we arrive easily at 
\begin{equation*}\J (a(Q))=(p-q)\sum \limits_{I=(i_1, \ldots, i_p)}^{} e_{i_1}^{\flat } \wedge \ldots \wedge e_{i_p}^{\flat}  \wedge (\bJ Q)(e_{i_1}, \ldots, e_{i_p}).
\end{equation*}
Applying $\J$ once more time while going through the same steps yields $\J^2 a(Q)=-(p-q)^2a(Q)$ and the proof is completed.
\end{proof}The main technical observation in this section is 
\begin{pro} \label{p21}
The following hold: 
\begin{itemize}
\item[(i)]
The total alternation map $a : \lambda^pV \ho \lambda^qV \to \Lambda^{p+q}V$ is injective for any $p \neq q$;
\item[(ii)]
The kernel of $a : \lambda^pV \otimes \lambda^qV \to \Lambda^{p+q}V$ is contained in $\lambda^pV \otimes_{2} \lambda^qV$.
\end{itemize}
\end{pro}
\begin{proof}
(i) If $Q$ belongs to  $\lambda^pV \ho \lambda^qV$ and $X$ is in $V$ we define $Q_X$ and $Q^X$ in $\lambda^{p-1}V \ho \lambda^{q}V$ and $\lambda^{p}V \ho \lambda^{q-1}V$ respectively by 
\begin{equation*}Q_X=Q(X, \cdot) \ \ \mbox{and} \ \ Q^X=X \lrcorner Q.
\end{equation*}
It is easy to see that those are well defined. Assume now that $a(Q)=0$.  Then 
\begin{align*}0=X \lrcorner a(Q)=&\sum_{i_1,\ldots ,i_p}X\lrcorner (e_{i_1}^{\flat}\wedge \ldots \wedge e_{i_p}^{\flat})\wedge Q(e_{i_1},...,e_{i_p})\\
&+(-1)^p\sum_{i_1,\ldots,i_p}e_{i_1}^{\flat}\wedge \ldots \wedge e_{i_p}^{\flat} \wedge (X\lrcorner Q(e_{i_1},..., e_{i_p}))\\
&=p\sum_{i_1,\ldots,i_{p-1}}e_{i_1}^{\flat}\wedge \ldots \wedge e_{i_{p-1}}^{\flat} \wedge Q(X,e_{i_1},..., e_{i_{p-1}})\\
&+(-1)^p\sum_{i_1,\ldots,i_p}e_{i_1}^{\flat}\wedge \ldots \wedge e_{i_p}^{\flat} \wedge Q^X(e_{i_1},..., e_{i_{p}})\\
&=p a(Q_X)+(-1)^p a(Q^X).
\end{align*}
By the previous Lemma $a(Q_X)$ is in $\lambda^{p-1,q}V$ whilst $a(Q^X)$ belongs to $\lambda^{p,q-1}V$ hence both must vanish since elements of distinct spaces as $p \neq q$. Now 
an induction argument leads directly to the proof of the claim in (i). \\
To prove (ii) we first note that $\lambda^pV \otimes \lambda^qV=(\lambda^pV \otimes_1 \lambda^qV) \oplus (\lambda^pV \otimes_2 \lambda^qV)$
and the claim follows from Lemma \ref{l31}.
\end{proof}
\vspace{-2.7mm}
Let $L : \Lambda^{\star}V \to \Lambda^{\star}V$ be the exterior multiplication with the K\"ahler form $\omega=g(J \cdot, \cdot)$. Recall that the space $\Lambda^{\star}_0V$ of primitive forms is defined to be the 
kernel of $L^{\star}$, the adjoint of $L$ w.r.t. the inner product $g$. We consider the operators 
$P_k : \Lambda^rV \times\Lambda^sV\rightarrow\Lambda^{r+s-2k}V$ defined by
$$P_k(\alpha,\beta):=\sum \limits_{1 \le i_1 \ldots i_k \le 2n}^{}(e_{i_1}\lrcorner \ldots e_{i_k}\lrcorner\alpha)\wedge(Je_{i_1}\lrcorner \ldots \lrcorner Je_{i_k}\lrcorner\beta)$$
for all $(\alpha, \beta)$ in $\Lambda^rV \times\Lambda^sV$ and where $\{e_i, 1 \le i \le 2n\}$ is some orthonormal basis in $V$. Clearly, $P_0(\alpha, \beta)=\alpha \wedge \beta$ for all $(\alpha, \beta)$ in $\Lambda^rV \times\Lambda^sV$ and 
moreover
\begin{pro}\label{pk} For any $\alpha\in\Lambda^rV$ and $\beta\in\Lambda^sV$, we  have
\begin{itemize}
\item[(i)]  $L^{\star}P_k(\alpha,\beta)=P_k(L^{\star}\alpha,\beta)+P_k(\alpha,L^{\star}\beta)+(-1)^{r-k-1}P_{k+1}(\alpha,\beta)$ for all $k \geq 0$;
\item[(ii)] $(L^{\star})^p(\alpha \wedge \beta)=(-1)^{\frac{p(p-1)}{2}} p!\langle \alpha, J\beta\rangle $ for any primitive $p$-forms $\alpha$ and $\beta$.
\end{itemize}
\end{pro}
\begin{proof}(i) Let $\alpha \in \Lambda^rV$ and $\beta \in \Lambda^sV$. Then
\begin{align*}  L^{\star}P_k(\alpha,\beta)&=\frac{1}{2}\sum_{i,i_1 \ldots i_k}Je_i\lrcorner  e_i\lrcorner ((e_{i_1}\lrcorner \ldots e_{i_k}\lrcorner\alpha)\wedge(Je_{i_1}\lrcorner \ldots Je_{i_k}\lrcorner\beta))\\
&=\frac{1}{2}\sum_{i,i_1 \ldots i_k}Je_i\lrcorner((e_i\lrcorner e_{i_1}\lrcorner \ldots e_{i_k}\lrcorner\alpha)\wedge(Je_{i_1}\lrcorner \ldots Je_{i_k}\lrcorner\beta))\\
&+\frac{1}{2}(-1)^{r-k}\sum_{i,i_1 \ldots i_k}Je_i\lrcorner((e_{i_1}\lrcorner \ldots e_{i_k}\lrcorner\alpha)\wedge(e_i\lrcorner Je_{i_1}\lrcorner \ldots Je_{i_k}\lrcorner\beta))\\
&=P_k(L^{\star}\alpha,\beta)+\frac{1}{2}(-1)^{r-k-1}\sum_{i_1 \ldots i_{k+1}}(e_{i_1}\lrcorner \ldots e_{i_{k+1}}\lrcorner\alpha)\wedge(Je_{i_1}\lrcorner \ldots Je_{i_{k+1}}\lrcorner\beta)\\
&+\frac{1}{2}(-1)^{r-k}\sum_{i_1 \dots i_{k+1}}(Je_{i_1}\lrcorner e_{i_2}\lrcorner \ldots e_{i_{k+1}}\lrcorner \alpha)\wedge(e_{i_1}\lrcorner Je_{i_2}\lrcorner \ldots Je_{i_{k+1}}\lrcorner\beta)\\
&+P_k(\alpha,L^{\star}\beta)
\end{align*}
\ndt and the claim in (i) follows. \\
To prove (ii) we first obtain by induction from (i) that $(L^{\star})^p(\alpha \wedge \beta)=(-1)^{\frac{p(p-1)}{2}}P_p(\alpha, \beta)$ whenever $\alpha, \beta$ belong to $\Lambda^p_0V$. To conclude 
it is enough to directly use the definition of $P_p$ to get $P_p(\alpha, \beta)=p!\langle \alpha, J\beta\rangle $.
\end{proof}
\subsection{Formal Sasakian metrics}
Part of the algebraic facts developed above can be also used to describe completely the cohomology algebra of a geometrically formal, Sasakian metric. For an introduction 
to Sasakian geometry, the odd dimensional analogue of K\"ahler geometry, we refer the reader to \cite{BoGa}.
\begin{teo} \label{Sasformal1} Let $(M^{2n+1},g)$ be a compact Sasakian manifold. If the metric $g$ is formal then $b_{p}(M)=0$ for all $1 \leqslant p \leqslant 2n $, in other words $M$ is a real cohomology sphere.
\end{teo}
\begin{proof}
Recall that the tangent bundle of $M$ splits as $TM=\V \oplus H$ an orthogonal direct sum where $\V$ is spanned by the so-called Reeb vector field, to be denoted by $\z$. The contact 
distribution $H$ admits a $g$-compatible complex structure $J: H \to H$ which moreover satisfies $d\theta=\omega$ where $\theta$ is the $1$-form dual to $\z$ and 
$\omega=g(J \cdot, \cdot)$. We call a differential $p$-form horizontal, and denote the corresponding space by $\Lambda^pH$ if the interior product with $\z$ vanishes.  Now let 
$d_H : \Lambda^{\star}H \to \Lambda^{\star}H$ be the projection of the usual exterior derivative $d$ onto $H$. If $d_H^{\star}$ is its formal adjoint w.r.t. to the restriction of $g$ 
on $H$, we have  (see \cite{ta}) on $\Lambda^{p}M=\Lambda^{p}H \oplus \biggl [\theta \wedge \Lambda^{p-1}H \biggr ]$
\begin{equation} \label{block}
d^{\star}=\biggl ( \begin{array}{cc}
d_H^{\star} & -{\mathcal{L}}_{\z}\\
L^{\star} & -d_H^{\star}
\end{array} \biggr ) 
\end{equation}
where ${\mathcal{L}}_{\z}$ denotes the Lie derivative. As a last reminder, we mention that the extension of $J$ to $\Lambda^{\star}H$ defined as in the previous section 
preserves the space of harmonic forms. 
$\\$
Let now $\alpha$ be a harmonic form on $M$. It is a known fact that if $0 \leqslant p \leqslant n$, every harmonic form $\alpha$ on $M$ is horizontal and invariant by the 
Reeb vector field. Moreover, $\alpha$ must be primitive, that is $L^{\star} \alpha=0$. Using the formality assumption on $g$ we obtain 
that $\alpha \wedge J \alpha$ is still harmonic. Since this is a horizontal form, invariant under the Reeb vector field it follows from \eqref{block} that 
$L^{\star}(\alpha \wedge J \alpha)=0$. We conclude that $\alpha$ vanishes by means of Proposition \ref{pk}, (ii).
\end{proof}
The proof of Theorem \ref{Sasformal} in the introduction is now complete.
\section{Holomorphic forms with harmonic squares}
Let $(M^{2n},g,J)$ be a compact K\"ahler manifold and consider a harmonic $p$-form $\Omega$ in $\lambda^pM$, that is of type $(0,p)+(p,0)$. It is a well known fact, see \cite{goldberg} for 
instance, that $\Omega$ must be holomorphic, that is 
\begin{equation} \label{h1}\nabla_{JX}\Omega=\nabla_X (\bJ \Omega)
\end{equation}for all $X$ in $TM$. Together with $\Omega$ comes 
$S : \Lambda^{p-1}M \to \Lambda^1M$ defined by $S(X_1,...,X_{p-1})=\Omega(X_1,...,X_{p-1}, \cdot)$. That $\Omega$ has real type $(0,p)+(p,0)$ translates into
\begin{equation} \label{h2}(S(JX_1,...,X_{p-1}))^{\sharp}=-J(S(X_1,...,X_{p-1}))^{\sharp}
\end{equation}
whenever $X_1,...,X_{p-1}$ belong to $TM$ and where for any $1$-form $\theta$, $\theta^{\sharp}$ denotes the associated vector field with respect to the metric $g$. Let now 
$Q : \Lambda^{p-1}M \to \lambda^pM$ be given by
\begin{equation*}Q(X_1,...,X_{p-1})=\nabla_{(S(X_1,\ldots,X_{p-1}))^{\sharp}}\Omega
\end{equation*}for all $X_1,...,X_{p-1}$ in $TM$. The next Lemma provides information about the complex type of $Q$.
\begin{lema} \label{l1h}The tensor $Q$ belongs to $\lambda^{p-1}M \ho \lambda^{p}M$.
\end{lema}
\begin{proof}
Follows immediately from (\ref{h1}) and (\ref{h2}).
\end{proof}
\begin{pro}\label{paraomega}Let $\Omega$ in $\lambda^pM$ be a harmonic form. If the metric $g$ is formal, then
\begin{equation}\nabla_{(S(X_1,\ldots,X_{p-1}))^{\sharp}}\Omega=0
\end{equation}holds, for all $X_1,\ldots,X_{p-1}$ in $TM$. 
\end{pro}
\begin{proof} Let $\{e_i, 1 \le i \le 2n \}$ be a geodesic frame at a point $m$ in $M$. If $p$ is even $\Omega \wedge \Omega$ is harmonic and we have at $m$ 
\begin{equation*}
\begin{split}0=&-d^{\star}(\Omega \wedge \Omega)=\sum \limits_{i=1}^{2n} e_i \lrcorner \nabla_{e_i}(\Omega \wedge \Omega) \\
=&2 \sum \limits_{i=1}^{2n} e_i \lrcorner ( \nabla_{e_i}\Omega \wedge \Omega)=2 \sum \limits_{i=1}^{2n} \nabla_{e_i}\Omega \wedge (e_i \lrcorner \Omega)
\end{split}
\end{equation*}since $\Omega$ is itself co-closed. In other words $a(Q)=0$ and we conclude by means of Lemma \ref{l1h} and Proposition \ref{p21} that $Q=0$. If $p$ is odd the harmonicity of 
$\Omega \wedge \bJ \Omega $ gives 
\begin{equation*}
\begin{split}
0=-d^{\star}(\Omega \wedge \bJ \Omega)&=\sum \limits_{i=1}^{2n} e_i \lrcorner (\nabla_{e_i}\Omega \wedge \bJ \Omega+\Omega \wedge \nabla_{e_i}\bJ \Omega) \\
& =\sum \limits_{i=1}^{2n}  -\nabla_{e_i}\Omega \wedge (e_i \lrcorner \bJ \Omega)+(e_i \lrcorner \Omega) \wedge \nabla_{e_i}(\bJ \Omega)
\end{split}
\end{equation*}
where we took into account the co-closedeness of $\Omega$ and $\bJ \Omega$. Now $\nabla_{e_i} \bJ \Omega =\nabla_{Je_i} \Omega$ hence
\begin{equation*}
\begin{split}
0&=\sum \limits_{i=1}^{2n} -\nabla_{e_i}\Omega \wedge (Je_i \lrcorner \Omega)+(e_i \lrcorner \Omega)\wedge \nabla_{Je_i}\Omega\\
&=-2\sum \limits_{i=1}^{2n}\nabla_{e_i}\Omega \wedge (Je_i \lrcorner \Omega).
\end{split}
\end{equation*}
This is easily reinterpreted to say that $a(\mathbb{J}Q)=0$ and then Lemma \ref{l1h} together with Proposition \ref{p21} leads to the vanishing of $Q$ and hence to the claimed result.
\end{proof}
\begin{rema} \label{lesshyp}From the proof of the result above we see that it actually holds for harmonic 
forms $\Omega$ in $\lambda^pM$ such that 
$\Omega \wedge \Omega$ ($p$ even) resp. $\Omega \wedge \bJ \Omega$ ($p$ odd) are co-closed.
\end{rema}
We need now to recall some facts about the algebraic structure of harmonic forms of type $(1,1)$. 
\begin{pro} \cite{nagy} \label{kahler}
Let $(M^{2n},g,J)$ be a compact  K\"ahler manifold such that the metric $g$ is formal. If $\alpha=g(F \cdot, \cdot)$ is harmonic in $\lambda^{1,1}M$ then we have 
an orthogonal and $J$-invariant splitting 
\begin{equation*}
TM=\bigoplus \limits_{i=0}^{p} E_i
\end{equation*}
which is preserved by $F$ and such that $F=\lambda_i J_i$ on $E_i$, for all $0 \leqslant i \leqslant p$. Here $J_i$ are almost complex structures on $E_i$ and $\lambda_i$ are real 
constants, for $0 \leqslant i \leqslant p$. 
\end{pro}
Now we would like to conclude from Proposition \ref{paraomega} that $\Omega $ is actually parallel. This is eventually seen to be the case 
if $\Omega $ is non-degenerate at every point of the manifold. To rule out the general case we must  study the null distribution of $\Omega$. For each 
$m$ in $M$ define $\V_m=\{X \in T_mM : X \lrcorner \Omega=0\}$. Our first concern is to show that $m \to \V_m$ gives a smooth, {\it{constant}} rank distribution on $M$.
\begin{lema} \label{nullrk} 
Let $(M^{2n},g,J)$ be a compact K\"ahler manifold such that the metric $g$ is formal. If $\Omega$ in $\lambda^pM$ is harmonic the following hold
\begin{itemize}
\item[(i)] the distribution $\V$ is of constant rank;
\item[(ii)] both distributions $\V$ and $H=\V^{\perp}$ are integrable and $H$ is totally geodesic.
\end{itemize}
\end{lema}
\begin{proof}
(i) Let $\alpha_{\Omega}$ in $\lambda^{1,1}M$ be defined by $\alpha_{\Omega}(X,Y)=\langle JX \lrcorner \Omega, Y \lrcorner \Omega\rangle $ for all $X,Y$ in $TM$. Because $g$ is formal 
we have that $(L^{\star})^{p-1}(\Omega \wedge J \Omega)$ is a harmonic two form. On the other hand side, from Proposition \ref{pk}, (i) it follows by induction that $(L^{\star})^{p-1}(\Omega \wedge J \Omega)=(-1)^{\frac{(p-2)(p-3)}{2}}
P_{p-1}(\Omega,J\Omega)$ by also using that $\Omega$ is primitive. Now a direct computation using the definition of $P_{p-1}$ shows that 
\begin{equation*}
\begin{split}
P_{p-1}(\Omega,J\Omega)(X,Y)=&(-1)^{p-1}(p-1)!(\langle X \lrcorner \Omega, JY \lrcorner \Omega \rangle-\langle Y \lrcorner \Omega, JX \lrcorner \Omega \rangle)\\
=&2(-1)^p(p-1)! \alpha_{\Omega}(X,Y)
\end{split}
\end{equation*}
for all $X,Y$ in $TM$. We conclude that $\alpha_{\Omega}$ is a harmonic form of type $(1,1)$ hence 
the formality of $g$ and Proposition \ref{kahler} ensure that $\alpha_{\Omega}$ has constant rank. By a positivity argument the nullity of $\alpha_{\Omega}$ coincides with that of $\Omega$ and 
the claim is proved.\\
(ii) $\V$ (hence $H$) is $J$-invariant since $\alpha_{\Omega}$ lives in $\lambda^{1,1}M$. By (i) we obtain a globally defined splitting $TM=\V \oplus H$ which is therefore orthogonal and $J$-invariant. From 
the definition of $\V$ it follows by an orthogonality argument that the distribution $H$ is spanned by $S(X_1,...,X_{p-1})$ with $X_1,...,X_{p-1}$ in $TM$ hence 
\begin{equation} \label{pomega}
\nabla_X \Omega=0 \ \mbox{for all} \ X \in H
\end{equation}
by Proposition \ref{paraomega}.
Taking now a 
direction, say $V$ in $\V$ gives that $\nabla_XV$ belongs to $\V$ and this shows the total geodesicity hence the integrability of $H$. The integrability of $\V$ is an 
easy consequence of the closedeness of $\Omega$. Indeed, taking $X_1,\ldots,X_{p-1}$ in $H$ and $V,W$ in $\V$, we have 
\begin{equation*}
\begin{split}
0=d\Omega(X_1,...,X_{p-1}, V,W)&=\sum \limits_{i=1}^{p-1}(-1)^{i+1}(\nabla_{X_i}\Omega)(X_1,...,\widehat{X_i},\ldots,X_{p-1},V,W)\\
&-(\nabla_V \Omega)(X_1,...,X_{p-1},W)+(\nabla_W \Omega)(X_1,\ldots,X_{p-1},V)\\
&=\Omega(X_1,...,X_{p-1},[V,W]).
\end{split}
\end{equation*}
Since $\Omega$ vanishes on $\V$ by the definition of the latter it follows that $[V,W] \lrcorner \Omega=0$ and our integrability claim follows by using again the definition of $\V$.
\end{proof}
To prove the parallelism of $\Omega$, which amounts to having $\V$ totally geodesic we need to establish one more fact. Recall \cite{To} that the transversal Ricci tensor 
$Ric^H:H \to H$ of the totally geodesic distribution $H$ is defined by
\begin{equation*}
g(Ric^HX,Y)=\sum \limits_{i}^{} R(X,e_i,Y,e_i)
\end{equation*}
for all $X,Y$ in $H$ and local orthonormal frames $\{e_i\}$ in $H$. When $\V$ integrates to give a Riemannian submersion, which is always true locally, $Ric^H$ corresponds to the usual Ricci tensor of the base manifold.
\begin{lema}\label{RicHvanish}The transversal Ricci tensor $\Ric^H$ of the distribution $H$ vanishes.
\end{lema}
\begin{proof} For any $\alpha$ in $\Lambda^2M$ and for all $\varphi$ in $\Lambda^{\star}M$ let us define 
\begin{equation*}
[\alpha,\varphi]=\sum \limits_{i=1}^{2n} e_i \lrcorner \alpha  \wedge e_i \lrcorner \varphi
\end{equation*}
where $\{e_i, 1 \le i \le 2n\}$ is some local orthonormal frame in $TM$. 
Since $H$ is totally geodesic, after differentiation of \eqref{pomega} in directions coming from $H$ we get $[R(X,Y),\Omega]=0$ for all $X,Y$ in $H$. Since $V \lrcorner \Omega=0$ for $V$ in $\V$ it follows that 
$\sum \limits_{i}R(X,Y)e_i \wedge e_i \lrcorner \Omega=0$ for all $X$ in $H$ and where $\{e_i\}$ is a local orthonormal frame in $H$, to be fixed in what follows. Therefore we get 
\begin{equation*}
\begin{split}
0=&\sum \limits_{j,i}^{} e_j \lrcorner (R(X,e_j)e_i \wedge e_i \lrcorner \Omega)=Ric^HX \lrcorner \Omega-\sum \limits_{j,i} R(X,e_j)e_i \wedge e_j \lrcorner e_i \lrcorner \Omega\\
=&Ric^HX \lrcorner \Omega+\frac{1}{2}\sum \limits_{j,i} R(e_j,e_i)X \wedge e_j \lrcorner e_i \lrcorner \Omega
\end{split}
\end{equation*}
for all $X$ in $H$, where for obtaining the second line we used the algebraic Bianchi identity for $R$. As consequences of the K\"ahler condition and of the fact that $\Omega$ is in $\lambda^pM$ we have that 
$R(JX,JY)=R(X,Y)$, whilst $JX \lrcorner JY \lrcorner \Omega=-X \lrcorner Y \lrcorner \Omega$ for all $X,Y$ in $TM$. Hence the last sum above vanishes and we end up with $Ric^HX \lrcorner \Omega=0$ for all 
$X$ in $TM$ whence the claim, since $\Omega$ is non-degenerate on $H$.
\end{proof}
At the same time, the situation when $Ric^H$ vanishes is well described by the following
\begin{teo}\label{oldth} \cite{nagy} Let $(M^{2n},g,J)$ be a compact K\"ahler manifold equipped with a Riemannian foliation with complex leaves. If the the foliation is 
transversally totally geodesic with nonnegative transversal Ricci tensor then it has to be totally geodesic, therefore locally a Riemannian product.
\end{teo}$\\$
{\bf{Proof of Theorem \ref{holo}}} Since $Ric^H$ vanishes, it follows by Theorem \ref{oldth} that $\V$ 
is totally geodesic, hence parallel w.r.t. the Levi-Civita connection $\nabla$. This implies immediately the parallelism of $\Omega$, by means of \eqref{pomega}. The local product decomposition of $(M^{2n},g,J)$ follows by using 
the deRham splitting theorem for the $\nabla$-parallel decomposition $TM=\V \oplus H$, combined with Lemma \ref{RicHvanish}.
\section{Harmonic $2$-forms}
We shall develop in this section the general Riemannian counterpart of Proposition \ref{kahler}. From now on, we shall use the metric to identify a 
$2$-form $\alpha$ with a skew-symmetric endomorphism $A$ of $TM$; explicitly $\alpha=g(A \cdot, \cdot)$. Moreover, the space $\A$ is the space of 
skew-symmetric endomorphisms of $TM$ which are associated to an element of $\H^2(M,g)$. If $\varphi$ belongs to $\Lambda^{\star}M$ let $L_{\varphi} : \Lambda^{\star}M 
\to \Lambda^{\star}M$ be given as exterior multiplication by $\varphi$ and let $L^{\star}_{\varphi}$ be the adjoint of $L_{\varphi}$.
\begin{pro}\label{stab} Let $M^n$ be geometrically 
formal and let $g$ be a formal metric on $M$. We have :
$$ A_2A_1A_3+A_3A_1A_2 \in \A$$ whenever $A_i, 1 \leqslant i \leqslant 3$ belong to $\A$. 
\end{pro}
\begin{proof}Let $\alpha$ belong to $\H^2(M,g)$. 
Since $g$ is formal and $L_{\alpha}^{\star}$ is up to sign equal to $\star L_{\alpha} \star$ it follows that both $L_{\alpha}$  and 
$L_{\alpha}^{\star}$ preserve the space of harmonic forms of $(M,g)$. Therefore, if $\alpha_i, 1 \leqslant i \leqslant 3$ belong to $\H^2(M,g)$ then 
$ L_{\alpha_1}^{\star}L_{\alpha_2}\alpha_3$ is an element of $\H^2(M,g)$. Let $A_i, 1 \leqslant i \leqslant 3$ be the skew-symmetric endomorphisms associated to the forms 
$\alpha_i, 1 \leqslant i \leqslant 3$ and let 
$\{e_i, 1 \leqslant i \leqslant n \}$ be a local orthonormal basis in $TM$. We shall now compute 
$$L_{\alpha_1}^{\star}L_{\alpha_2}\alpha_3=\frac{1}{2}
\sum \limits_{i,j=1}^{n} \alpha_1(e_i,e_j) e_j \lrcorner \biggl [ e_i \lrcorner (\alpha_2 \wedge \alpha_3 )\biggr ]$$
But 
$$e_j \lrcorner \biggl [ e_i \lrcorner (\alpha_2 \wedge \alpha_3 ) \biggr ]=
\alpha_2(e_i,e_j)\alpha_3-(e_i \lrcorner \alpha_2) \wedge (e_j \lrcorner \alpha_3)+(e_j \lrcorner \alpha_2) \wedge
 (e_i \lrcorner \alpha_3)+\alpha_3(e_i,e_j) \alpha_2.$$
Further computation yields, after some elementary manipulations 
$$L_{\alpha_1}^{\star}L_{\alpha_2}\alpha_3=\langle\alpha_1, \alpha_2\rangle  \alpha_3+
\langle\alpha_1, \alpha_3\rangle  \alpha_2 +\langle A_3A_1A_2+A_2A_1A_3 \cdot, \cdot\rangle $$
\end{proof}
In what follows we shall say that a symplectic form on $M$ is compatible with the metric $g$ if its associated skew-symmetric endomorphism defines an almost complex structure on $M$.
\begin{pro}\label{decomp} Let $M^n$ be geometrically formal and let $g$ denote a formal metric on $M$. Moreover, let $\alpha$ belong to $\H^2(M,g)$ with associated endomorphism $A$ in 
$\A$. Then:
\begin{itemize} 
\item[(i)] The eigenvalues of $A^2$ are constant with eigenbundles of constant rank;
\item[(ii)] Let $\mu_i$ be (the pairwise distinct) 
eigenvalues of $A^2$, with $\mu_0=0$ and let $E_i$ be the eigenbundles of $A^2$ corresponding to $\mu_i$. Then for $1\leqslant i\leqslant p$,  
$E_i$ is of even dimension and we have an orthogonal decomposition $$\alpha=\sum_{i=1}^p {\sqrt{-\mu_i}} \omega_i$$\ndt where $\omega_i, 1 \le i \le p$ belong to 
$\H^2(M,g)$. Moreover, $\omega_i=g(J_i \cdot, \cdot)$ on $E_i$, for some $g$-compatible almost complex structure $J_i$ on $E_i, 1 \leqslant i\leqslant p$;
\item[(iii)] if $\alpha$ is non-degenerate then $g$ admits a compatible symplectic form.
\end{itemize}
\end{pro}
\begin{proof}
(i) From Proposition \ref{stab} we get by induction that $A^{2k+1}$ belongs to $\A$ whenever $A$ is in $\A$. Since $\A$ is finite dimensional, there exists 
$P\in\mathbb{R}[X]$ so that $P(A^2)=0$ and moreover by using the symmetry of $A^2$ the polynomial $P$ can be supposed to have only real and simple roots $\mu_i, 1 \le i \le p$. 
Let $m_i$ be the dimension of the $\mu_i$-eigenbundle, $1 \le i \le p$. To see that $m_i, \mu_i$, 
$0\leqslant i\leqslant p$ are constant over $M$, we use the fact that $A^{2k+1}$ belongs to $\A$ for any $k\in\mathbb{N}$ by Proposition \ref{stab} and from the fact that elements in $\A$ have 
pointwisely constant scalar products we deduce that $Tr (A^{2k})=-\langle A^{2k-1} , A\rangle=c_k$ for some constant $c_k$ and for any integer $k$. It follows that 
$\sum \limits_{i=1}^{p}m_i\mu_i^k=c_k$ for all $k$ in $\mathbb{N}$ hence
this Vandermonde system leads to the constancy of the functions $m_i,\mu_i$, $1\leqslant i\leqslant p$.\\
(ii) With the notation $\lambda_i=\sqrt{-\mu_i}$, the orthogonal projection of $\alpha$ on $E_i$ is given by $\lambda_i\omega_i$ where $\omega_i=g(J_i \cdot, \cdot)$ for some 
almost complex structure  $J_i$ on $E_i$, $1\leqslant i\leqslant p$. Now
\begin{equation*}
g(A^{2k+1}\cdot, \cdot)=\displaystyle\sum_{i=1}^p \lambda_i^{2k+1}\omega_i
\end{equation*} 
is harmonic for all natural $k$ and by an argument similar to the one used in the proof of the Proposition 3.1 of \cite{nagy} we 
deduce that $\omega_i$ belong to $\H^2(M,g)$.\\
(iii) By (ii) the form $\sum \limits_{i=1}^p \omega_i$ belongs to $\H^2(M,g)$ and it is $g$-compatible if $\alpha$ is non-degenerate.
\end{proof}
The technical advantage of Proposition \ref{decomp} is essentially to say that all distributions appearing as ranges or kernels of harmonic $2$-forms are of constant rank over the manifold, and in this respect 
they can -as we shall see in the next section-be treated as algebraic objects. 
\subsection{$6$-dimensions}
We shall present here a geometric application of the algebraic facts from the previous section. More precisely, we are going to obtain  
sufficient conditions for a geometrically formal $6$-manifold to admit a compatible symplectic structure. We need first to make a number of preliminary results.
\begin{lema} \label{hforms}
Let $M^n$ be geometrically formal and let $g$ be a formal metric on $M$. Let $\alpha$ belong to $\H^2(M,g)$ with kernel $\V$ and such that on $H={\V}^{\perp}, \alpha=g(J \cdot, \cdot)$ for some 
almost complex structure $J$ of $H$. Then for any $\phi$ in $\H^{p}(M,g)$ we have that $\phi^{ij}$ belongs to $\H^{p}(M,g)$ where for any $i,j$ with $i+j=p$ we have 
denoted by $\phi^{ij}$ the orthogonal projection of $\phi$ onto $\Lambda^i\V \widehat{\otimes} \Lambda^jH \subseteq \Lambda^pM$. Here $\Lambda^i\V \widehat{\otimes} \Lambda^jH$ is the image of 
$\Lambda^i\V \otimes \Lambda^jH$ in $\Lambda^{i+j}M$ under the antisymmetrisation map.
\end{lema}
\begin{proof}
We first note that
\begin{equation*}
L^{\star}_{\alpha}(\psi \wedge \alpha)=\frac{1}{2}(-1)^p (dim\  H )\psi+(L^{\star}_{\alpha}\psi) \wedge\alpha+(-1)^p Q \psi
\end{equation*} whenever $\psi$ is a $p$-form on $M$, where the operator $Q$ is given by $Q \psi=\sum \limits_{e_i \in H}^{}(e_i \lrcorner \psi) \wedge e^i$ for 
an arbitrary local frame $\{e_i\}$ in $H$. Hence $Q$ preserves the space of harmonic forms and on the other hand a standard computation shows that the nonzero eigenvalues 
of $Q$ on $\Lambda^pM$ are $(-1)^{p-1} j$ for $1 \leqslant j\leqslant dim\ H$ and $i=p-j\leqslant dim\ \V$ with corresponding eigenbundles $\Lambda^i\V \widehat{\otimes} \Lambda^jH$.
But formality actually implies that all powers of $Q$ preserve $\H^p(M,g)$, and the claim follows.
\end{proof}
\begin{lema} \label{notsymp}
Let $M^6$ be geometrically formal and let $g$ be a formal metric on $M$. If $g$ does not admit a compatible symplectic form then every non-zero harmonic $2$-form on $M$ has $4$-dimensional 
kernel. 
\end{lema}
\begin{proof} Let $\alpha \neq 0$ belong to $\H^2(M,g)$. It cannot be non-degenerate for Proposition \ref{decomp}, (iii)  would imply the existence of a $g$-compatible symplectic form.
It remains to see that $\alpha$ cannot have $2$-dimensional kernel. Arguing by contradiction, let us suppose that $\V=Ker(\alpha)$ 
is $2$-dimensional, so that $H=\V^{\perp}$ is of dimension $4$. Moreover, from $\alpha$ we get again by using  Proposition \ref{decomp} a harmonic $2$-form 
$\alpha^{\prime}=g(J \cdot, \cdot)$ on $H$ for some almost complex structure $J$ on $H$. Then $\alpha^{\prime}+\star(\alpha^ {\prime} \wedge \alpha^{\prime})$ gives a 
globally defined symplectic form on $M$, compatible with $g$, hence the desired contradiction. 
\end{proof}
In what follows the distribution spanned by an orthonormal system of vector fields $\{X_1, \ldots,X_q\}$ on $M$ shall be denoted by $(X_1,\ldots,X_q)$.
\begin{pro}\label{classtart} Let $M^6$ be geometrically formal with $b_1(M)=0$ and $b_2(M) \geqslant 2$. If $g$ is a formal metric on $M$ which does not admit a compatible symplectic form we must 
have $b_2(M)=2, b_3(M)=6$.
\end{pro}
\begin{proof}
Let $\alpha \neq 0$ belong to $\H^2(M,g)$. By Lemma \ref{notsymp} the distribution $\V=Ker(\alpha)$ must be $4$-dimensional, so after constant rescaling 
$\alpha$ can be written as $\alpha=g(J \cdot, \cdot)$ where $J$ is an almost complex structure on the plane distribution $H=\V^{\perp}$. We now note there are no 
non-zero harmonic $2$-forms contained in $\Lambda^2\V$, for by Lemma \ref{notsymp} any such form must have $4$-dimensional kernel and hence must vanish.
It follows then from Lemma \ref{hforms} that $\H^2(M,g)$ is contained in $(\Lambda^1\V \widehat{\otimes} \Lambda^1H) \oplus \mathbb{R}\alpha$. Further on, because $b_2(M) \geqslant 2$, there 
must be a non-zero $\beta$ in $\Lambda^1\V \widehat{\otimes} \Lambda^1H$, and again by Lemma \ref{notsymp} this has $4$-dimensional kernel to be denoted by $\V^{\prime}$.  By rescaling 
if necessary we may also assume that $\beta$ is of unit length. \\
Let now $F_1$ and $F_2$ be the orthogonal projections of $H^{\prime}=({\V}^{\prime})^{\perp}$ onto $\V$ and $H$ respectively. $F_1$ is not the zero space because otherwise we would have 
$H^{\prime} \subseteq H$ hence $\beta$ in $\Lambda^2H$, an absurdity. We cannot have $F_2=\{0\}$ neither: it would imply that $H^{\prime} \subseteq \V$ hence
$\beta\in\Lambda^2\V$ which is again impossible. Therefore, both of $F_1$ and $F_2$ have rank at least $1$ and given that $H^{\prime}=F_1 \oplus F_2$ and $H^{\prime}$ has 
rank $2$, their respective ranks must actually equal $1$. Since the manifold is 
oriented, every real line bundle over $M$ is trivial and this leads to the existence 
of a globally defined orthonormal frame $\{\z, e_2\}$ on $H^{\prime}$, spanning $F_1$ and $F_2$.  Since $\beta$ belongs to $\Lambda^2H^{\prime}$, it follows that 
\begin{equation*}
\beta=e^2 \wedge \z.
\end{equation*}
Now the orthogonal complement of $(e_2)$ in $H$ is $1$-dimensional, hence trivial as a real line bundle. Therefore it is spanned by some a unit vector field, say $e_1$, and since 
$\alpha$ belongs to $\Lambda^2H$ we get  
\begin{equation*}
\alpha=e^1 \wedge e^2.
\end{equation*}
Pick now a non-zero harmonic $3$-form $T$ on $M$. By Lemma \ref{hforms} the components $T^{11}$ in $\Lambda^3\V$ and $T^{12}=\theta\wedge\alpha$ in 
$\Lambda^1\V \widehat{\otimes} \Lambda^2H$ of $T$ are harmonic. But $\star L_{\alpha} T^{11}$ and $L^{\star}_{\alpha}T^{12}=\theta$ are harmonic $1$-forms and since $b_1(M)=0$ these 
$1$-forms are vanishing fact which implies the nullity of $T^{11}$ and $T^{12}$. Hence $T$ can be written as 
\begin{equation*}
T=\omega_1 \wedge e^1+\omega_2 \wedge e^2
\end{equation*}
with $\omega_k, k=1,2$ in $\Lambda^2\V$. Again, $L_{\phi}T$ resp. $L^{\star}_{\phi}T$ vanish for any harmonic $2$-form $\phi$ because $b_1(M)=0$, hence from 
$L_{\beta}T=0$ and $L^{\star}_{\beta}T=0$ we get that 
\begin{equation*}
\z \wedge \omega_1=0, \ \z \lrcorner \omega_2=0.
\end{equation*}
It follows easily that harmonic $3$-forms on $M$ are contained in a rank $6$ sub-bundle of $\Lambda^3M$, thus using that scalar products of harmonic $3$-forms are (pointwisely) constant we obtain that $b_3(M) \leqslant 6$. 
Since $M$ has nowhere vanishing vector fields, it has vanishing Euler characteristic, and from $b_1(M)=0, b_2(M) \geqslant 2$ we get
\begin{equation*}
b_3(M)=2(1+b_2(M)) \geqslant 6
\end{equation*}
showing that actually $b_2(M)=2$ and $b_3(M)=6$.
\end{proof}
\begin{teo} \label{intermed}Let $M^6$ be geometrically formal with $b_1(M) \neq 1$ and $b_2(M) \geqslant 2$. If $g$ is a formal metric on $M$ which does not admit a compatible symplectic form then either: \\
(i) $M$ has the real cohomology algebra of $\mathbb{T}^3 \times S^3$ \\
or \\
(ii) $b_1(M)=0, b_2(M)=2, b_3(M)=6$.
\end{teo}
\begin{proof}
In view of the Proposition above it suffices to treat the cases when $b_1(M) \neq 0$. Again, we do a case by case discussion. Let $\V$ be the distribution spanned by the harmonic 
$1$-forms and let $\z_k, 1 \leqslant k \leqslant b_1(M)$ be a frame of harmonic $1$-forms in $\V$. As an immediate consequence of 
Lemma \ref{hforms} and of the fact that $H={\V}^{\perp}$ does not contain, by definition, harmonic $1$-forms it follows that harmonic $2$-forms are contained in $\Lambda^{2}\V \oplus 
\Lambda^2H$. \\ 
If $b_1(M)=2$, $H$ is of rank $4$ and since $b_2(M) \geqslant 2$ there must be a non-zero harmonic $2$-form contained in $\Lambda^2H$. In view of Lemma \ref{notsymp} it has rank $4$ kernel and therefore vanishes, a contradiction. 
$\\$ 
Suppose now that $b_1(M)=3$ so that $H$ is of rank $3$. If $\alpha$ is a non-zero harmonic $2$-form contained in $\Lambda^2H$, then $\z_1 \lrcorner \z_2 \lrcorner \z_3 \lrcorner (\star \alpha)$ is a non-zero harmonic form in 
$\Lambda^1H$ which is a contradiction. Therefore $\H^2(M,g) \subseteq \Lambda^2\V$ and similarly, by using Lemma \ref{hforms} we get 
$\H^3(M,g) \subseteq \Lambda^{3}\V \oplus \Lambda^3H$. It is now straightforward that $M$ has the cohomology algebra of $\mathbb{T}^3 \times S^3$.\\
If $b_1(M)=4$, then $\z_1 \wedge \z_2+\z_3 \wedge \z_4+\star(\z_1 \wedge \z_2 \wedge \z_3 \wedge \z_4)$ is a compatible symplectic form, a contradiction.\\
Now we cannot have $b_1(M)=5$ (\cite{Kot}) and when $b_1(M)=6$ there exists an orthonormal frame of harmonic $1$-forms hence a compatible symplectic structure, a contradiction. This finishes the proof of 
the Theorem.
\end{proof}
The proof of Theorem \ref{symp}, when $b_2(M) \geqslant 3$ follows now immediately from the above. 
\begin{rema}The proof of Proposition \ref{classtart} can also be adapted to show that if $g$ is a formal metric on $M^6$ which does not admit 
a compatible symplectic structure then $b_3(M) \leqslant 6$ when $b_1(M)=0, b_2(M)=1$. 
\end{rema}
\subsection{The case when $b_1=0, b_2=2, b_3=6$}
We shall examine now the case when the geometrically formal manifold $M^6$ has a formal metric $g$ which does not admit a compatible symplectic form and moreover $b_1(M)=0, b_2(M)=2, b_3(M)=6$.  We have seen that harmonic 
$2$-forms must be of the form $e^{12}=e^1 \w e^2, e^2 \w \z$ for some orthonormal system $e_1,e_2,\z$ in $TM$. Let us denote by $E$ the rank $3$ distribution orthogonal to $e_1,e_2, \z$. It inherits 
a transversal volume form, i.e a nowhere vanishing $3$-form $\nu_E$ in $\Lambda^3E$ given by $\nu_E=\star (e^{12} \w \z)$. We shall write $\star_E : \Lambda^{\star}E 
\to \Lambda^{\star}E$ for the Hodge star operator obtained when $E$ is equipped with the restriction of the metric $g$ and orientation given by $\nu_E$.
\begin{lema} \label{streq}The following hold : 
\begin{equation*}
\begin{split}
de^1=& A \w e^1+B \w e^2+\lambda e^{12}\\
de^2=&q\z \w e^1-A \w e^2+\mu e^{12}\\
d\z=&A \w \z-\mu e^1 \w \z+e^2 \w D
\end{split}
\end{equation*}
where $A,B,D$ are $1$-forms on $E \oplus ( \z)$ and $\lambda, q, \mu$ are functions on $M$.
\end{lema}
\begin{proof}
Because $e^{12}$ is closed we get $de^1 \w e^2=de^2 \w e^1$ and it follows that none of $de^1, de^2$ can have components in $\Lambda^{2}(e_1,e_2)^{\perp}$.  Therefore one can write 
\begin{equation*}
\begin{split}
de^1=& A \w e^1+B \w e^2+\lambda e^{12}\\
de^2=&C \w e^1+D^{\prime} \w e^2+\mu e^{12}\\
\end{split}
\end{equation*}
for some one-forms $A,B,C,D^{\prime}$ in $\Lambda^1(e_1,e_2)^{\perp}$ and some smooth functions $\lambda,\mu$ on $M$. Now the remaining information contained in $de^1 \w e^2=de^2 \w e^1$ is 
that $D^{\prime}=-A$. Since $e^2 \w \z$ is equally closed we have $de^2\w \z=d \z \w e^2$ hence $de^2 \w \z \w e^2=0 $ leading to $C \w \z=0$. Thus we may write $C=q\z$ for some smooth function 
$q$ on $M$. Moreover, by an argument already used for $e^{12}$, $d\z$ has no component in $\Lambda^2 (e_2,\z)^{\perp}$ hence after a small computation we can fully rewrite the closedeness 
of $e^2 \w \z$ as 
$$ d\z=A \w \z-\mu e^1 \w \z+e^2 \w D+\nu e^{12}$$
for some one form $D$ on $E \oplus (\z)$ and a smooth function $\nu$ on $M$. Now the harmonicity of $e^{12}$ tells us that 
$$ 0=d^{\star}e^{12}=d^{\star}e^1 \cdot e^2-[e_1,e_2]-d^{\star}e^2 \cdot e_2$$
in other words the distribution $(e_1,e_2)$ is integrable. Henceforth, $\nu=d\z(e_1,e_2)=-<\z, [e_1,e_2]>$ vanishes and our Lemma is proved.
\end{proof}
\begin{coro} \label{corint}
\begin{itemize}
\item[(i)] The distribution $E$ is integrable.
\item[(ii)] The distributions $(e_1,e_2)$ and $(e_2,\z)$ are integrable as well.
\end{itemize}
\end{coro}
\begin{proof}
(i) By inspecting the structure equations in the Lemma above, we see that either of $d\z, de^1, de^2$ vanish on $\Lambda^2E$ and the claim follows.\\
(ii) follows by arguments similar to the last part of the proof of the Lemma \ref{streq}.
\end{proof}
We shall now bring into consideration the fact that $b_3(M)=6$. Let 
\begin{equation} \label{h3f}
T_1, T_2, T_3 , \star T_1, \star T_2, \star T_3 
\end{equation}
be an (pointwisely) orthonormal basis in $\mathcal{H}^3(M,g)$. From the proof of Proposition \ref{classtart} we must have 
\begin{equation*}
T_k=(e^1 \w \z) \w \alpha_k+e^2 \w \star_E \beta_k
\end{equation*}
where $\alpha_k, \beta_k$ belong to $\Lambda^1E$ for all $1 \le k \le 3$. The next Lemma recasts the orthogonality of the system 
\eqref{h3f} into a simpler algebraic form. 
\begin{lema} \label{algi3} 
For $1 \le k \le 3$ we define $\gamma_k=\alpha_k+i\beta_k$ in $\Lambda^1(E, \mathbb{C})$. We have 
\begin{equation*}
\begin{split}
\star_E\gamma_1&= k \overline{\gamma}_2 \w \overline{\gamma}_3\\
\star_E \gamma_2&=-k \og_1 \w \og_3 \\
\star_E \gamma_3&=k \og_1 \w \og_2
\end{split}
\end{equation*} 
for some smooth function $k : M \to \mathbb{C}$ such that $\vert k \vert=1$ and $k \og_1 \w \og_2 \w \og_3=\nu_E$.
\end{lema}
\begin{proof} 
The Hodge star operator of the forms $T_k, 1 \le k \le 3$ reads 
\begin{equation*}
\star T_k=-(e^1 \w \z) \w \beta_k+e^2 \w \star_E \alpha_k
\end{equation*} 
and the orthonormality of \eqref{h3f} is equivalent with the following 
\begin{equation*} 
\begin{split}
&\vert \alpha_k\vert^2+\vert \beta_k \vert^2=1\\
& <\alpha_i, \alpha_j>+<\beta_i, \beta_j>=0, i \neq j \\
&  <\alpha_i, \beta_j>=<\alpha_j, \beta_i>
\end{split}
\end{equation*}
It is easy to see that $\{\gamma_i,1 \le i \le 3\}$ gives a basis of $\Lambda^1(E, \mathbb{C})$ (not orthonormal though) and then $\{\gamma_i \w \gamma_j : 1 \le i \neq j \le 3\}$ is a basis 
in $\Lambda^2(E, \mathbb{C})$. Of course, by using complex conjugation we obtain another set of basis in the above mentioned spaces. We now compute 
\begin{equation*}
\begin{split}
\star_E \gamma_j \w \og_j=&(\star_E \alpha_j +i\star_E \beta_j) \w (\alpha_j -i\star_E \beta_j)\\
=& (\star_E\alpha_j \w \alpha_j+\star_E \beta_j)+i(\star_E \beta_j \w \alpha_j-\star_E \alpha_j \w \beta_j)\\
=&\nu_E.
\end{split}
\end{equation*}
Very similarly, we also find that $\star_E \gamma_j \w \og_p=0$ for $p \neq j$ and the result follows. That $\vert k \vert=1$ follows routineously by taking norms.
\end{proof}
The triple of $1$-forms  $(\gamma_1, \gamma_2, \gamma_3)$ has also an internal symmetry, of particular relevance for what follows. Write 
$\gamma=\biggl ( \begin{array}{c}
\gamma_1 \\
\gamma_2 \\
\gamma_3
\end{array}\biggr )$ and then notice the transition formula $\gamma=P \overline{\gamma}$ for some $P=(P_{ij}, 1 \le i, j \le 3) : M \to M_3(\mathbb{C})$. This 
is possible because both $\gamma$ and $\overline{\gamma}$ give basis in $\Lambda^1(E, \mathbb{C})$. It follows immediately that $P \overline{P}=I_3$ holds and moreover 
from the definition of $P$ we see that it is symmetric, i.e. $P=P^{T}$. To exploit the closedeness the frame \eqref{h3f} we need the following preliminary 
\begin{lema} \label{extd}If $\alpha$ belongs to $\Lambda^{\star}E$ we have 
\begin{equation*}
d\alpha=d_E\alpha+\z \w L^{E}_{\z} \alpha+e^1 \w (L^{E}_{e_1} \alpha+\z \w R \lrcorner \alpha)+e^2 \w L^{E}_{e_2} \alpha
\end{equation*}
where $d_E$ denotes the orthogonal projection of $d$ onto $\Lambda^{\star}E$ and for any vector field $X$ in $E, L_X^{E}$ is the orthogonal projection of 
the Lie derivative $L_X\alpha$ onto $\Lambda^{\star}E$. Moreover, the vector field $R$ in $E$ is given by the projection on $E$ of $[e_1,\z]$.  
\end{lema}
\begin{proof} Follows eventually by expanding $d$ along the decomposition 
$$\Lambda^{\star}M=\Lambda^{\star}E \otimes \Lambda^{\star}(e_1,e_2,\z)$$
while making use of the integrability of the distributions listed in Corollary \ref{corint}.
\end{proof}
Let us denote by $\hat{A}, \hat{B}, \hat{D}$ the components on $E$ of the $1$-forms $A,B,D$, so that $A=\hat{A}+x\z, B=\hat{B}+y\z, D=\hat{D}+z\z$ for some smooth functions 
$x,y,z$ on $M$.
\begin{lema}\label{sysh} The harmonicity of the forms $T_k, 1 \le k \le 3$ is equivalent with the following system of equations: 
\begin{itemize}
\item[(i)] $d_E \gamma_k=-2\hat{A} \wedge \gamma_k-iq \star_E \gamma_k$
\item[(ii)] $d_E (\star_E \gamma_k )=\hat{A} \w \star_E \gamma_k $ 
\item [(iii)] $L_{\z}^E ( \star_E \gamma_k)-x \star_E \gamma_k-i\hat{B} \w \gamma_k=0$
\item [(iv)] $L_{e_1}^E ( \star_E \gamma_k )+\mu \star_E \gamma_k-i\hat{D} \w \gamma_k=0$
\item [(v)] $L_{e_2}^E \gamma_k+(z-\lambda)\gamma_k-iR \lrcorner \star_E \gamma_k=0$
\end{itemize}
for $1 \le k \le 3$.
\end{lema}
\begin{proof}
For any $1 \le k \le 3$ the closedeness of the forms $T_k$ is equivalent with 
\begin{equation*}
\begin{split}
0=dT_k=& d(e^1 \w \z) \w \alpha_k+e^1 \w \z \w d\alpha_k \\
+& de^2 \w [\star_E \beta_k]-e^2 \w d[\star_E \beta_k].
\end{split}
\end{equation*}
Using now Lemma \ref{extd} we obtain further 
\begin{equation*}
\begin{split}
0&=d(e ^1 \w \z) \w \alpha_k+de^2 \w \star_E \beta_k\\
&+ e^1 \w \z \w \biggl [ d_E \alpha_k+e^2 \w L_{e_2}^E \alpha_k\biggr ] \\
&-e^2 \w d_E (\star_E \beta_k)-e^2 \w \z \w L_{\z}^E(\star_E \beta_k)+e^{12} \w L_{e_1}^E(\star_E \beta_k)\\
&-e^{12} \w \z \w (R \lrcorner \star_E \beta_k).
\end{split}
\end{equation*}
But accordingly to Lemma \ref{streq} we eventually get 
\begin{equation*}
d(e^1 \w \z)=2\hat{A} \w e^1 \w \z+\hat{B} \w e^2 \w \z-\hat{D} \w e^{12}+(\lambda-z)e^{12} \w \z  
\end{equation*}
hence after identifying the components of $e^1 \w \z, e^2 \w \z, e^{12}, e^{12} \w \z, e^2$ we find the system of equations 
\begin{equation*}
\begin{split}
& 2\hat{A} \w \alpha_k-q \star_E \beta_k+d_E \alpha_k=0\\
&\hat{B} \w \alpha_k+x \star_E \beta_k-L_{\z}^E(\star_E \beta_k)=0\\
&-\hat{D} \w \alpha_k+\mu \star_{E}\beta_k+L_{e_1}^E (\star_E \beta_k)=0\\
& (\lambda-z)\alpha_k-L^E_{e_2} \alpha_k -R \lrcorner \star_E \beta_k=0\\
& \hat{A} \wedge \star_E \beta_k=d_E(\star_E \beta_k)
\end{split}
\end{equation*}
But the forms $\star T_k, 1 \le k \le 3$ are closed as well, in other words the system above has the symmetry $(\alpha_k,\beta_k) \to (\beta_k, -\alpha_k) $. It is now straightforward to rephrase these 
by means of the complex valued forms $\gamma_k , 1 \le k \le 3$.
\end{proof}
We are now in position to examine the geometric consequences imposed by our initial situation. 
\begin{lema} The following hold: 
\begin{itemize}
\item[(i)] $\hat{A}=0$;
\item[(ii)] $d_Ek=0$.
\end{itemize}
\end{lema}
\begin{proof}
We will prove both claims at the same time. Using Lemma \ref{sysh}, (i) we compute 
\begin{equation*}
\begin{split}
d_E(\gamma_2 \w \gamma_3)=&-4\hat{A} \w \gamma_2 \w \gamma_3-iq (\star_E \gamma_2 \w \gamma_3-
\star_E \gamma_3 \w \gamma_2)\\
=&-4 \hat{A} \w \gamma_2 \w \gamma_3
\end{split}
\end{equation*}
by using standard properties of the Hodge star operator. But from (ii) of the same Lemma, actualised by Lemma \ref{algi3} one infers that 
$$ d_E(k\og_2 \w \og_3)=k\hat{A} \w \og_2 \w \og_3. $$
It follows that $(5\hat{A}-k^{-1}d_Ek) \w \og_2 \w \og_3=0$ and repeating the procedure for the other two equations in Lemma \ref{sysh}, (i) we arrive easily to 
$5\hat{A}-k^{-1}d_Ek=0$. But $\hat{A}$ is real valued whilst $k^{-1}d_Ek$ belongs to $\Lambda^1(E, i\mathbb{R})$ since $\vert k \vert=1$ and the proof of the Lemma follows.
\end{proof}
We examine the rest of the equations in Lemma \ref{sysh}. For a triple $\alpha=\biggl ( \begin{array}{c}
\alpha_1\\
\alpha_2 \\
\alpha_3
\end{array}\biggr )$ of one forms in $\Lambda^1(E, \mathbb{C})$ we consider the triple of $2$-forms in $\Lambda^2(E, \mathbb{C})$ given by $\alpha \times \alpha=
\biggl( \begin{array}{c}
\alpha_2 \w \alpha_3 \\
\alpha_3 \w \alpha_1 \\
\alpha_1 \w \alpha_2 
\end{array} \biggr )$. Note that in the new notation Lemma \ref{algi3} now reads 
\begin{equation} \label{vprod1}
\star_E \gamma=k \og \times \og
\end{equation} and after taking the conjugate 
we also get 
\begin{equation} \label{vprod2}
\star_E \og=k^{-1}\gamma \times \gamma
\end{equation}
since $\overline{k}=k^{-1}$. For any $\alpha=\sum \limits_{k=1}^3 \alpha_k \gamma_k$ in 
$\Lambda^1(E, \mathbb{C})$ we consider the matrix 
\begin{equation*}
r_{\alpha}=\biggl ( \begin{array}{ccc}
0 & \alpha_3 & -\alpha_2 \\
-\alpha_3 & 0 & \alpha_1 \\
\alpha_2 & -\alpha_1 & 0 
\end{array} \biggr )
\end{equation*}
Note that $r_{\alpha}^{T}=-r_{\alpha}$ and we shall let $r_{\alpha}$ operate on triple of forms in $\Lambda^k(E,\mathbb{C}), k=1,2$ by matrix multiplication. Moreover, a straightforward computation shows that 
$\alpha \w \gamma=\biggl ( \begin{array}{c}
\alpha \w \gamma_1 \\
\alpha \w \gamma_2 \\
\alpha \w \gamma_3 
\end{array} \biggr )=r_{\alpha}(\gamma \times \gamma)$. These observations allow now to bring the remaining equations into final form. 
\begin{lema} \label{reint}
The following hold 
\begin{itemize}
\item[(i)] $L_{\z}^E (\star_E \gamma)-x \star_E \gamma-ik r_{\hat{B}}(\star_{E} \og)=0$
\item[(ii)] $L_{e_1}^E (\star_E \gamma)+\mu \star_E \gamma-ik r_{\hat{D}} (\star_{E} \og)=0$ 
\item[(iii)] $L_{e_2}^E \gamma+(z-\lambda)\gamma+ikr_{\eta} \og=0$\\
where the $1$-form $\eta$ in $\Lambda^1E$ is given as $\eta=g(R, \cdot)$.
\end{itemize}
\end{lema}
\begin{proof}
We shall prove only (i) the other two claims being entirely analogous. Indeed, writing (iii) of Lemma \ref{sysh} in matrix form we have 
\begin{equation*}
L_{\z}^E (\star_E \gamma)-x \star_E \gamma-i\hat{B} \w \gamma=0.
\end{equation*}
But $\hat{B} \w \gamma=r_{\hat{B}} (\gamma \times \gamma)=kr_{\hat{B}} (\star_{E} \og)$ by \eqref{vprod2} and we are done.
\end{proof}
\begin{pro}\label{obstr} The following hold: 
\begin{itemize}
\item[(i)] $L_{e_1}P=L_{e_2}P=L_{\z}P=0$
\item[(ii)] $P\overline{r_{\hat{B}}}P+k^2r_{\hat{B}}=0$
\item[(iii)] $P\overline{r_{\hat{D}}}P+k^2r_{\hat{D}}=0$
\item [(iv)] $P\overline{r_{\eta}}P+k^2r_{\eta}=0$.
\end{itemize}
\end{pro}
\begin{proof}
Taking the conjugate in (i) of Lemma \ref{reint} we get 
\begin{equation} \label{istep1}
L^E_{\z} (\star_E \og)-x \star_E \og+ik^{-1} \overline{r_{\hat{B}}}(\star_E \gamma)=0.
\end{equation}
Now $\star_E \gamma=\star_E(P \og)=P(\star_E \og)$ hence (i) of Lemma \ref{reint} gives 
\begin{equation*}
(L_{\z}^E P) \star_E \og+PL_{\z}^E(\star_E \og)+xP(\star_E \og)+ik r_{\hat{B}}(\star_{E} \og)=0
\end{equation*}
Substituting here the expression of $L_{\z}^E(\star_E \og)$ as given by \eqref{istep1} we obtain further 
\begin{equation*}
\begin{split}
(L_{\z}^E P) \star_E \og +& P \biggl [ x \star_{E} \overline{\gamma}-ik^{-1}\overline{r}_{\hat{B}}(\star_E \gamma)\biggr ]\\
-& (xP+ik r_{\hat{B}}) \star_{E} \og=0
\end{split}
\end{equation*}
whence
\begin{equation*} 
(L_{\z}^E P -ik^{-1}P\overline{r_{\hat{B}}}P-ik r_{\hat{B}} ) \star_E \og=0
\end{equation*}
where we have used once more that $\gamma=P\og$. Given that $\star_E \og$ gives a basis in $\Lambda^2(E, \mathbb{C})$ we infer that 
$$ L_{\z}^E P -ik^{-1}P\overline{r_{\hat{B}}}P-ik r_{\hat{B}} =0.
$$
But $P$ is symmetric and $r_{\hat{B}}$ is skew-symmetric therefore $P\overline{r_{\hat{B}}}P$ is skew-symmetric as well, hence identifying the symmetric resp. the skew-symmetric 
part in the equation above we arrive at $L_{\z}^E P=0$ and $P\overline{r_{\hat{B}}}P+k^2r_{\hat{B}}=0$. The other two claims in (i) and assertions in (iii) and (iv) are proved by applying a 
completely similar procedure to the equations in (ii) and (iii) of Lemma \ref{reint}.
\end{proof}
\begin{coro} \label{vanishi2}We must have $\hat{B}=\hat{D}=\eta=0$.
\end{coro}
\begin{proof}
We first  work out the equation in (ii) of Lemma \ref{obstr}. It implies that 
$$ (P\overline{r_{\hat{B}}}P) \star_E \og+k^2 r_{\hat{B}}(\star_E \og)=0.
$$
Now $r_{\hat{B}}(\star_E \og)=k^{-1}r_{\hat{B}}(\gamma \times \gamma)=k^{-1} \hat{B} \w \gamma$. On the other hand we have 
\begin{equation*}
\begin{split}
(P\overline{r_{\hat{B}}}P) \star_E \og=&(P\overline{r_{\hat{B}}})\star_E (P\og)\\
=&P\overline{r_{\hat{B}}}\star_E \gamma=kP\overline{r_{\hat{B}}} (\og \times \og )\\
=&kP(\overline{\hat{B}} \w \og)\\
=&k \hat{B} \w P\og=k \hat{B} \w \gamma
\end{split}
\end{equation*}
since $\hat{B}$ is real valued. Altogether $(k+k^{-1}k^2) \hat{B} \w \gamma=0$ whence the vanishing of $\hat{B}$ since $\vert k \vert=1$. The vanishing of $\hat{D}$ resp. $\eta$ follows now from 
(iii) resp. (iv) of Lemma \ref{obstr} by using the same argument.
\end{proof}
We now continue the study of the distribution $(e_1,e_2, \z)$. 
\begin{lema} \label{codiff}The following hold:
\begin{itemize}
\item[(i)] $d^{\star}e^1=-\mu$;
\item[(ii)] $d^{\star} e^2=\lambda=-z$;
\item[(iii)] $d^{\star}\z=x$.
\end{itemize}
\end{lema}
\begin{proof}
First of all we update Lemma \ref{streq} to 
\begin{equation} \label{strlast}
\begin{split}
de^1=& x\z \w e^1+y \z \w e^2+\lambda e^{12}\\
de^2=&q\z \w e^1-x\z \w e^2+\mu e^{12}\\
d\z=&-\mu e^1 \w \z+ze^2 \w \z
\end{split}
\end{equation}
by using that $\hat{A}=\hat{B}=\hat{D}=0$. \\
(i) since $e^{12}$ is harmonic we have 
$$ 0=d^{\star}(e^{12})=d^{\star}e^1\cdot e_2-[e_1,e_2]-d^{\star}e^2 \cdot e_1
$$
hence $d^{\star}e^1=<[e_1,e_2],e_2>=-de^2(e_1,e_2)=-\mu$ and $d^{\star}e^2=-<[e_1,e_2],e_1>=de^1(e_1,e_2)=\lambda$. This proves (i) and the first half of (ii) To prove the rest 
it is enough to repeat the argument above starting from $d^{\star}(e^2 \w \z)=0$.
\end{proof}
\begin{teo} \label{b026} A geometrically formal manifold $M^6$ with $b_1(M)=0, b_2(M)=2, b_3(M)=6$ and formal metric $g$ must admit a $g$-compatible symplectic structure.
\end{teo}
\begin{proof}Suppose that there is no $g$-compatible symplectic structure on $M$. Then our whole previous discussion applies and based upon it we will obtain a contradiction.\\ 
We proceed first towards updating the expressions of the Lie derivatives of $\gamma, \star_E \gamma$ as given by Lemma \ref{reint}. Since $k^2=det(P)$ and $P$ has no 
Lie derivatives in the direction of $(e_1,e_2, \z)$ it follows that $L_{e_1}k=L_{e_2}k=L_{\z}k=0$. Therefore, (i) of Lemma \ref{reint} gives 
$$ L_{\z}^E(\og \times \og)-x \og \times \og=0.$$
Note that actually $L^E_{\z}\gamma=L_{\z}\gamma$ since $\eta$ (hence $R$) vanishes. A short computation using only that $\gamma$ gives a basis in $\Lambda^{1}(E, \mathbb{C})$ leads to 
\begin{equation*}
L_{\z} \gamma-\frac{x}{2}\gamma=0.
\end{equation*}
It follows that $L_{\z} (\gamma_1 \w \gamma_2 \w \gamma_3)=\frac{3x}{2}\gamma_1 \w \gamma_2 \w \gamma_3$ whence $L_{\z}\nu_E=\frac{3x}{2}\nu_E$. But $L_{\z}(e^{12} \w \z)=0$ as well, 
because $e^{12}, e ^{12} \w \z$ are closed (the latter after a computation based on \eqref{strlast}) and we get that the volume form 
$\nu_M=e^{12}\w \z \w \nu_E$ satisfies $L_{\z}\nu_M=\frac{3x}{2}\nu_M$. But 
$$ L_{\z}\nu_M=d(\z \lrcorner \nu_M)=-d^{\star}\z \cdot \nu_M=-x\nu_M$$
by Lemma \ref{codiff}, (iii) and it follows that we must have $x=0$. When working out, in the same spirit, the equation contained in (ii) of Lemma \ref{reint} we obtain that $\mu=0$. \\
Now (iii) of Lemma \ref{reint} ensures, as before, that $L_{e_2}\nu_E+3(z-\lambda)\nu_E=0$. At the same time 
$$ L_{e_2}(e^{12} \w \z)=-d(e^1 \w \z)=(-\lambda+z)e^{12}\w \z
$$
after making use of \eqref{strlast}. It follows that $L_{e_2}\nu_M=-2(z-\lambda) \nu_M=4\lambda \nu_M$ as $z=-\lambda$ by Lemma \ref{codiff}, (ii).  But once again from 
$L_{e_2}\nu_M=-d^{\star}e^2 \cdot \nu_M=-\lambda \cdot \nu_M$ we obtain that $\lambda=0$. \\
Inspecting now the structure equations in \eqref{strlast} we see that $d\z=0$ and again from Lemma \ref{codiff} $d^{\star}\z=0$, in other words $\z$ is a harmonic, nowhere 
vanishing $1$-form on $M$ which contradicts that $b_1(M)=0$.
\end{proof}
The proof of the Theorem \ref{symp} in the introduction is now complete.
\section{Formal metrics with maximal $b_2$}We study in this section geometrically formal manifolds $M^{n}$ having maximal second Betti 
number, i.e. $b_2(M)=\begin{pmatrix} n\\ 2\end{pmatrix}$. To prove Theorem  \ref{max2}, we split our discussion into two cases according to the parity of 
$n$.  
\begin{pro} \label{obsvmax}Let $M^n$ be geometrically formal and let $g$ be a formal metric on $M$. The following hold:
\begin{itemize}
\item[(i)] if $b_{p}(M)$ and $b_{q}(M)$ are maximal for $p+q\leqslant n$ then $b_{p+q}(M)$ is also maximal;
\item[(ii)]  if $b_{p}(M)$ and $b_{q}(M)$ are maximal for $0 \leq p \leq q \leq n$ and  then so is $b_{q-p}(M)$;
\item[(iii)]if $b_p(M)$ is maximal for some $ 1 \leqslant p \leqslant n-1$ and $(p,n)=1$ then $g$ is a flat metric. 
\end{itemize}
\end{pro} 
\begin{proof} 
(i) If $\{\alpha_i\}, \{\beta_j\}$ are $L^2$-orthonormal basis in $\H^p(M,g)$ and $\H^q(M,g)$ respectively then at each point of $M$ we obtain orthonormal basis in $\Lambda^pM$ and $\Lambda^qM$ respectively. 
It follows that $\Lambda^{p+q}M$ is spanned by forms of the type $\alpha_i \wedge \beta_j$ which are harmonic because the metric $g$ is formal. Since scalar products between harmonic forms are constant 
after Gramm-Schmidt orthonormalisation we obtain a basis in $\H^{p+q}(M,g)$.\\
(ii) By Hodge duality $b_{n-p}(M)$ is maximal hence by (i) so is $b_{n-p+q}(M)=b_{q-p}(M)$ whence the claim.\\
(iii) If $b_p(M)$ is maximal then for any integers $q$ and $k$, $1\leqslant k\leqslant n$ such that $pq\equiv k (\text{mod}\ \  n)$ $b_k(M)$ is also maximal by using (i). Since $(p,n)=1$ we arrive by means of (ii) at $b_1(M)$ maximal, and it follows
that $g$ is flat by Theorem \ref{thK}, (iii).
\end{proof}
Hence, when $n$ is odd and $b_2(M)$ is maximal $b_1(M)$ is maximal too and the metric $g$ is flat. Therefore we need only to consider the case when $n$ is even.
\subsection{Reduction to the symplectic case}
As an immediate consequence of Proposition \ref{decomp} we have :
\begin{pro} \label{redsymp} Let $M^{n}$ be a geometrically formal manifold with formal metric $g$ such that $b_2(M)$ is maximal and $n$ is even. Then $g$ admits a compatible almost K\"ahler structure, that is an almost complex 
structure $J$, which is compatible with $g$ and such that the $2$-form $g(J \cdot, \cdot)$ is closed.
\end{pro}
\begin{proof} We first claim that there 
exists a harmonic $2$-form $\alpha$ which is non-degenerate, that is $\alpha^k \neq 0, n=2k$ at some point $x$ of $M$. Indeed if $\varphi^k=0$ on $M$ for any $\varphi$ in $\H^2(M,g)$ then after polarisation we find 
$\varphi_1 \wedge \ldots \wedge \varphi_k=0$ whenever $\varphi_i, 1 \le i \le k$ belong to $\H^2(M,g)$. Since frames in $\H^2(M,g)$ give frames in the $\Lambda^2M$ it is easy to obtain a contradiction and the existence of $\alpha$ as above follows.
The claim is now proved by using (iii) in Proposition \ref{decomp}.
\end{proof}
\subsection{Proof of flatness}We consider hereafter a compact almost-K\"ahler manifold $(M^{n},g,J)$ ($n=2k$) such that $g$ is a formal metric and moreover $b_2(M)=\begin{pmatrix}n\\ 2\end{pmatrix}$. Let $\omega=g(J \cdot, \cdot)$ be the 
so-called K\"ahler form of the almost K\"ahler structure. We first remark that the bi-type splitting of $\Lambda^2M$ is preserved at the level of harmonic 
forms (note, by contrast with the K\"ahler case that this needs no longer be true in the case of an arbitrary almost K\"ahler manifold).
\begin{lema} \label{split2symp}Any harmonic $2$-form 
splits as $\alpha=\alpha_1+\alpha_2$ where the harmonic $\alpha_1, \alpha_2$ are in $\lambda^{1,1}M$ and $\lambda^2M$ respectively.
\end{lema}
\begin{proof}Pick $\alpha$ in $\Lambda^2M$, which splits as $\alpha=\alpha_1+\alpha_2$ with $\alpha_1$ in $\lambda^{1,1}M$ and 
$\alpha_2$ in $\lambda^2M$. Because of formality we can assume w.l.o.g. that $\alpha$ is primitive. Again the formality tells us that 
$L_{\alpha}^{\star}(\omega \wedge \omega)$ is harmonic and from the proof of Proposition \ref{stab} it follows that it is actually proportional to 
$\alpha_1-\alpha_2$.  This eventually proves the Lemma.
\end{proof}
Therefore, if $b_2(M)$ is maximal, both $\lambda^{1,1}M$ and $\lambda^2M$ are 
spanned by harmonic forms. We need now to see which geometric properties a harmonic $2$-form in $\lambda^2M$ must have. To do so, recall that 
the first canonical Hermitian connection $\bnabla$ of the almost K\"ahler $(g,J)$ is given by 
\begin{equation*}\bnabla_X=\nabla_X+\eta_X\end{equation*} for all $X$ in $TM$. Here $\nabla$ is the Levi-Civita connection of $g$ and 
$\eta_X=\frac{1}{2}(\nabla_XJ)J$ for all $X$ in $TM$ gives the intrinsic torsion 
of the $U(n)$-structure induced by $(g,J)$. The 
connection $\bnabla$ is metric and Hermitian, that is it preserves both the metric and the almost-complex structure. The almost K\"ahler condition i.e. that 
$d \omega=0$, when formulated in terms of the intrinsic torsion tensor $\eta$ reads  
\begin{equation} \label{ak}\langle\eta_XY,Z\rangle +\langle\eta_YZ,X\rangle +\langle\eta_ZX,Y\rangle =0
\end{equation}for all $X,Y,Z$ in $TM$. The latter also implies that 
$(g,J)$ is quasi-K\"ahler: 
\begin{equation} \label{qK}\eta_{JX}=\eta_X J
\end{equation}for all $X$ in $TM$. Moreover we have 
\begin{align}\label{commut} \eta_X J=-J\eta_X
\end{align}
in other words $\eta$ belongs to $\lambda^1M \otimes_{1} \lambda^2M$. The relations (\ref{ak}), (\ref{qK}) and (\ref{commut}) will be used implicitly in 
subsequent computations.
\begin{lema}\label{nablaf} Let $(M^{2k},g,J)$ be an almost-K\"ahler manifold and let $\alpha=g(F \cdot, \cdot)$ be harmonic in $\lambda^2M$. Then 
\begin{equation} \label{dprime}(\bnabla_{JX}F)JY+(\bnabla_XF)Y=-2\eta_{FX}Y
\end{equation}for all $X,Y$ in $TM$.
\end{lema}
\begin{proof}From $d \alpha=0$ we have that $a(\nabla \alpha)=0$. But $\nabla_X\alpha=\overline{\nabla}_X \alpha+\langle [F,\eta_X] \cdot, \cdot\rangle $ for all $X$ in 
$TM$ and moreover a 
simple computation based on (\ref{ak}) shows that  
\begin{equation*}a((X,Y,Z) \to \langle [F,\eta_X]Y,Z\rangle )=a((X,Y,Z) \to \langle \eta_{FX}Y,Z\rangle ).
\end{equation*}
Therefore $a(\overline{\nabla}\alpha+\eta_{F \cdot})=0$ and since the tensor under alternation belongs to $\lambda^{1}M \otimes \lambda^2M$ we use 
Proposition \ref{p21}, (ii) to conclude that it is actually in $\lambda^1M \otimes_{2} \lambda^2M$ and the proof of the claim follows by using the relations 
(\ref{qK}), (\ref{commut}).
\end{proof}
If $Q$ is an endomorphism of $M$, let us define the tensor $Q \bullet \eta$ by 
\begin{equation*}(Q \bullet \eta)(X,Y,Z)= \sigma_{X,Y,Z}\langle\eta_{QX}Y,Z\rangle 
\end{equation*}for all $X,Y,Z$ in $TM$, where $\sigma$ stands for the cyclic sum. Note that this is different from the usual action of $End(TM)$.
\begin{lema}Let $(M^{2k},g,J)$ be an almost-K\"ahler manifold and let $\alpha=g(F \cdot, \cdot)$ be harmonic in $\lambda^2M$ with harmonic square. Then 
\begin{equation} \label{algc}F^2 \bullet \eta=0.
\end{equation}
\end{lema}
\begin{proof}
That $d^{\star}(\alpha \wedge \alpha)=0$ translates after a calculation which parallels that in the proof of Proposition \ref{paraomega} into 
\begin{equation*}\sigma_{X,Y,Z} \langle(\nabla_{FX} F)Y,Z\rangle =0
\end{equation*} for all $X,Y,Z$ in $TM$. Rewritten by means of the canonical Hermitian connection and using (\ref{ak}) it yields
\begin{equation}\label{nablafcyc}
 \begin{split}&\langle(\bnabla_{FX}F)Y,Z\rangle +\langle(\bnabla_{FY}F)Z,X\rangle +\langle(\bnabla_{FZ}F)X,Y\rangle  \\
&+\langle\eta_XFY,FZ\rangle +\langle\eta_YFZ,FX\rangle +\langle\eta_ZFX,FY\rangle =0
\end{split}
\end{equation}
We shall exploit now the algebraic symmetries of the above. Changing $(Y,Z)$ in $(JY,JZ)$ and subtracting from the original equation implies 
\begin{equation*} 
\begin{split}
&2\langle(\bnabla_{FX}F)Y,Z\rangle -2\langle\eta_{X}FZ,FY\rangle \\
&+\langle(\bnabla_{FY}F)Z+(\bnabla_{JFY}F)JZ,X\rangle -\langle(\bnabla_{FZ}F)Y+(\bnabla_{JFZ}F)JY,X\rangle =0
\end{split}
\end{equation*}
or further, after using the relation \eqref{dprime}
\begin{equation} \label{derF}
\begin{split}
&\langle(\bnabla_{FX}F)Y,Z\rangle -\langle\eta_X FZ,FY\rangle \\
&-\langle\eta_{F^2Y}Z,X\rangle +\langle X, \eta_{F^{2}Z}Y\rangle =0.
\end{split}
\end{equation}
Now taking the cyclic sum and using \eqref{nablafcyc} we get the desired result.
\end{proof}
\begin{rema}
On an almost K\"ahler manifold $(M^{2k},g,J)$ a harmonic form $\alpha$ in $\lambda^2M$ with harmonic exterior powers needs not to be parallel w.r.t to the Levi-Civita connection of the metric 
$g$. This happens 
for instance when $\alpha=g(I \cdot, \cdot)$ for a $g$-compatible almost complex structure $I$ with $IJ+JI=0$,  which actually induces a complex-symplectic structure on $M$. Examples 
in this direction, which are not hyperk\"ahler, can be constructed on certain classes of nilmanifolds \cite{anna}.
\end{rema}
From the Lemma above we find by $J$-polarisation that 
$$ [F,G]\bullet \eta=0$$for all $F,G$ dual to harmonic forms in $\lambda^2M$. It is well known that the splitting $\mathfrak{so}(2k)=\mathfrak{u}(k) \oplus \mathfrak{m}$, where $\mathfrak{m}$ consists in elements of $\mathfrak{so}(2k)$ 
anti-commuting with $J$, is such that 
$[\mathfrak{m}, \mathfrak{m}]=\mathfrak{u}(k)$ for $k \geqslant 2$. Therefore,  if  $g$ is a formal metric on $M^{2k}$ and $b_2(M)$ 
is maximal, we get that $F \bullet \eta=0$ for all $F$ dual to forms in $\lambda^{1,1}M$ provided that $dim M \geqslant 6$.
\begin{lema} \label{van} 
If $\dim M\geqslant 6$, the intrinsic torsion $\eta$ must vanish identically.
\end{lema}
\begin{proof}
It is enough to prove the statement at an arbitrary point $m$ of $M$.
Pick an arbitrary unit vector $V$ in $T_mM$ and let $F$ be the skew-symmetric, $J$-invariant endomorphism of $TM$ which is $J$ on 
$E=\langle\{V,JV\}\rangle $ and vanishes on $H=E^{\perp}$. That $F \bullet \eta=0$ says
\begin{equation*}
\langle\eta_{FX}Y,Z\rangle +\langle\eta_{FY}Z,X\rangle +\langle\eta_{FZ}X,Y\rangle =0
\end{equation*}
for all $X,Y,Z$ in $TM$. It follows that $\langle\eta_V X,Y\rangle =0$ for all $X,Y$ in $H$, hence $\eta_V X$ is in $E$ for any $X\in H$. Moreover, since 
$\dim M \geqslant 6$, there exists a unit vector $U\in TM$ so that $(V,JV,U,JU,X,JX)$ is an orthogonal system. Let us consider the skew-symmetric, $J$-invariant endomorphism 
$G$ of $TM$ defined by $GV=U$, $GJV=JU$, $GU=-V$, $GJU=-JV$ and $G$ vanishes on $E'^{\perp}$ where $E'=\langle\{V,JV,U,JU\}\rangle $. Then 
\begin{equation*}\langle\eta_{GU}X,V\rangle +\langle\eta_{GX}V,U\rangle +\langle\eta_{GV}U,X\rangle =0
\end{equation*}\ndt 
This implies that 
$\langle\eta_V X,V\rangle=-\langle\eta_U X,U\rangle$.  Changing $V$ in $JV$ and using the $J$-anti-invariance of $\eta$ we get $\langle\eta_V X,V\rangle=0$. Then 
\begin{equation*}\label{nullite}
\eta_V X=0
\end{equation*}
for all $X\in H$ and $\eta_V X=\langle X,V\rangle \eta_V V+\langle X,JV\rangle \eta_V JV$ for all $X\in TM$. But from (\ref{nullite}) it follows that 
$\eta_V V=\eta_V JV=0$ and $\eta_V X=0$ for all $X\in TM$. 
\end{proof}
In other words $(g,J)$ is a K\"ahler structure and the flatness of the metric follows now from \cite{nagy}.  To finish the 
proof of Theorem \ref{max2} it remains to treat the case when $n=4$. In this situation, we notice that the bundles $\Lambda^{\pm}M$ of (anti) self-dual forms are 
trivialised by almost-K\"ahler structures satisfying the quaternionic identities and using the well-known Hitchin Lemma \cite{Hi} we obtain that 
$\Lambda^{\pm}M$ both contain a hyper-K\"ahler structure and this leads routineously to the flatness of the metric. \\
$\\$
{\bf{Acknowledgements:}} During the preparation of this paper, the research of P-A.N. was partly supported by the VW Foundation, through the program 
"Special geometries in mathematical physics", at the HU of Berlin, and an UoA grant. He is also 
grateful to the Institute \'E. Cartan in Nancy for warm hospitality during his visits there. We thank the referee and U. Semmelmann for useful suggestions on how to 
improve this work.


\begin{thebibliography}{99}
\bibitem{besse} A. L. Besse, \textit{Einstein Manifolds}, Springer Verlag, 1986.
\bibitem{deligne}
P. Deligne, Ph.Griffiths, J.Morgan, D.Sullivan,  \textit{Real homotopy type of K\"ahler manifolds}, Invent. Math. {\bf{29}} (1975), no.3, 245-274.
\bibitem{anna} 
A.Fino, H. Pedersen, Y. S. Poon and M. Weye Sorensen
\textit{Neutral Calabi-Yau structures on Kodaira Manifolds}, Commun. Math. Phys. 248 (2004) no. 2, 255-268.
\bibitem{BoGa}
Ch.P.Boyer, K.Galicki, \textit{Sasakian geometry}, Oxford University Press, 2008.
\bibitem{goldberg}
S.I.Godberg, \textit{Curvature and Homology}, Academic Press, 1962.
\bibitem{Kot} D.Kotschik, \textit{On products of harmonic forms}, Duke Math. J. {\bf 107}, no 3, (2001), 521-531.
\bibitem{fer1}
M. Fern\'andez, V.Mun\~oz, \textit{Formality of Donaldson submanifolds}, Math.Z. 
{\bf{250}} (2005), 149-175.
\bibitem{nagy}P-A.Nagy, \textit{On length and product harmonic forms in K\"ahler geometry},  Math. Z.  {\bf{254}} (2006),  no. 1, 199--218. 
\bibitem{ill}
J. Neisendorfer, T.J. Miller, \textit{Formal and coformal spaces}, Illinois J.Math.{\bf{22}} (1978), 565-580.
\bibitem{Hi} 
N.J.Hitchin, \textit{The self-duality equations on a Riemann surface}, Proc. London Math. Soc. {\bf{55}} (1987), 59-126.
\bibitem{To}
Ph. Tondeur, \textit{Geometry of foliations}, Birkh\"auser Verlag, 1997.
\bibitem{SMS}
S.M.Salamon, \textit{Riemannian geometry and holonomy groups}, Pitman Research 
Notes in Mathematics Series {\bf{201}}, 1989.
\bibitem{Sul}
D. Sullivan, \textit{Differential forms and the topology of manifolds} in {\it{Manifolds (Tokyo 1973)}}, ed. A. Hattori, Tokyo Univ. Press, 1975, 37-49.
\bibitem{W}B.Watson, \textit{Almost Hermitian submersions}, J.Diff.Geom. {\bf{11}} (1976), 147-165.
\bibitem{ta}
S. Tanaka, \textit{A differential geometric study on strongly pseudo-convex manifolds}, Lectures in Mathematics, Department of Mathematics, Kyoto University, No. 9. Kinokuniya Book-Store Co., Ltd., Tokyo, 1975. 
\bibitem{tanre}
D. Tanr\'e, \textit{Homotopie rationelle : Mod\`eles de Chen ,Quillen, Sullivan}, 
Lecture Notes in Math. {\bf{1025}}, Springer Verlag, 1983.
\end{thebibliography}
\end{document}